\documentclass[11pt,a4paper]{article}
\pdfoutput=1

\usepackage[utf8]{inputenc}
\usepackage[english]{babel}
\usepackage[left=2.5cm,right=2.5cm,top=2.5cm,bottom=3cm]{geometry}
\usepackage{amsmath}
\usepackage{amsfonts}
\usepackage{amssymb}
\usepackage{mathtools}

\usepackage{amsthm}
\usepackage{changepage}
\usepackage{chngcntr} 
\usepackage{faktor}
\usepackage{stmaryrd}
\usepackage[dvipsnames]{xcolor}

\usepackage{enumitem}
\usepackage{caption}
\usepackage{cancel}
\usepackage{easy-todo}
\newcommand{\say}[1]{``#1"}
\usepackage{import}

\usepackage{tikz}
\usepackage{tikz-cd}
\usetikzlibrary{shapes, calc, arrows, patterns, decorations.markings, decorations.pathreplacing, decorations.pathmorphing, decorations.shapes}
\pgfdeclarelayer{background} 
\pgfdeclarelayer{foreground}
\pgfsetlayers{background,main,foreground}

\usepackage{pifont} 
\newcommand{\cmark}{\ding{51}}
\newcommand{\xmark}{\ding{55}}

\newtheorem{thm}{Theorem}
\newtheorem*{thm*}{Theorem}
\newtheorem{lem}[thm]{Lemma}
\newtheorem{co}[thm]{Corollary}
\newtheorem{prop}[thm]{Proposition}

\theoremstyle{definition}
\newtheorem{defn}[thm]{Definition}

\newtheorem{rmk}[thm]{Remark}

\usepackage{makecell}

\newcommand{\N}{\mathbb{N}}		
\newcommand{\Z}{\mathbb{Z}}		
\newcommand{\R}{\mathbb{R}} 	

\newcommand{\ab}{\mathrm{ab}}
\newcommand{\la}{\left\langle}
\newcommand{\ra}{\right\rangle}

\newcommand{\ConvHull}{\mathrm{ConvHull}}

\newcommand{\CC}{\mathsf{CC}}

\newcommand{\CF}{\mathsf{CF}}
\newcommand{\Ac}{\mathcal A}
\newcommand{\Lc}{\mathcal L}
\DeclareMathOperator{\ev}{ev}
\newcommand{\rat}{\mathrm{rat}}
\newcommand{\alg}{\mathrm{alg}}

\newcommand{\standard}{\mathrm{std}}
\newcommand{\cubelike}{\mathrm{cube-like}}
\newcommand{\stndlike}{\mathrm{stnd-like}}

\renewcommand{\NG}{\N G}
\newcommand{\ZG}{\Z G}

\newcommand{\into}{\hookrightarrow}
\newcommand{\norm}[1]{\left\| #1 \right\|}
\newcommand{\abs}[1]{\left| #1 \right|}

\newcommand{\dif}{\mathrm d}
\DeclareMathOperator{\supp}{supp}
\newcommand{\num}{\mathrm{num}}

\usepackage{url}
\makeatletter
\newcommand{\address}[1]{\gdef\@address{#1}}
\newcommand{\email}[1]{\gdef\@email{\url{#1}}}
\newcommand{\emaill}[1]{\gdef\@emaill{\url{#1}}}
\newcommand{\@endstuff}{\par\vspace{\baselineskip}\noindent\small
	\begin{tabular}{@{}ll}\scshape\@address \\ \textit{E-mail address:} \@email \end{tabular}}
\AtEndDocument{\@endstuff}
\makeatother

\title{Dead ends and rationality of complete growth series}
\author{Pierre Alderic Bagnoud \and Corentin Bodart}
\address{Section de Mathématiques, Université de Genève, 1205 Genève, Switzerland}
\email{corentin.bodart@unige.ch}

\begin{document}
	\maketitle
	
	\begin{abstract}
	We are interested in the $\N G$-rationality and $\N G$-algebraicity of the complete growth series of finitely generated groups. It is shown that dead ends of arbitrarily large depths form an obstruction to $\N G$-rationality. In the case of the $3$-dimensional Heisenberg group $H_3(\Z)$, we prove that the complete series is not $\N G$-algebraic for any generating set. Dead ends are also used to show that complete growth series of higher Heisenberg groups are not $\N G$-rational for specific generating sets. Using a more general version of this obstruction, we prove that complete growth series of some lamplighter groups are not $\N G$-rational either. This work provides the first examples of groups exhibiting a difference between rationality of standard growth series, and rationality of complete growth series.
\end{abstract}

The idea of geometric group theory is to see a group $G$ as a geometric space: fix a finite symmetric generating set $S$, and endow $G$ with the word metric. Once this is done, we may be interested in the growth of the group, specifically the sequence $(a_n)_{n\ge 0}$ counting elements $g\in G$ at distance $\norm g_S=n$ of the identity. This sequence can be looked at under a \say{coarse} lens, and turned into a quasi-isometry invariant. In another direction, we may study finer properties of the sequence. This is usually done through the \emph{\say{numerical} growth series} of $G$
\[ G_\num(s) = \sum_{n\ge 0} a_n\cdot s^n \in \N[[s]]. \]
In this paper, we will be interested in an even richer sequence, namely $A_n=\sum_{\norm g=n} g$, seen as elements of the group semiring $\N G$. Just as in the standard case, the sequence $(A_n)_n$ gives rise to a growth series: the \emph{complete growth series} of $G$
\[G(s)=\sum_{n\ge 0} A_n\cdot s^n\in \N G[[s]]. \]
Note that both constructions are strongly related. Indeed, there exists an \emph{augmentation map} $\epsilon\colon\N G\to\N$, defined as the morphism sending all group elements to $1$. This morphism naturally extends into a morphism $\epsilon\colon\N G[[s]]\to \N[[s]]$. Observe that $\epsilon(A_n)=a_n$ hence $\epsilon(G(s))=G_\num(s)$.

\bigbreak

A series with coefficients in a semiring $R$ is \emph{$R$-rational} if it satisfies a certain system of linear equations (see \S\ref{sec:NG_rat_alg} for a proper definition). Similarly, the series is \emph{$R$-algebraic} if it satisfies a certain system of polynomial equations.

\bigbreak

The interest of rationality for numerical growth series is that their coefficients are relatively easy to compute: they satisfy a linear recurrence equation. Moreover precise asymptotics are known. Similarly, knowing that the complete growth series of a group is $\ZG$-rational says that the list of elements of length $n$ satisfies some linear recurrence relations. In some sense, being $\NG$-rational is even computationally better. (See Remark \ref{sec1:rmk_compu} for a discussion.)

\bigbreak

Proving complete growth series of groups are $\NG$-rational (and sometimes computing them) has been central since their introduction. Indeed, the two thesis of Liardet and Nagnibeda which introduced complete growth series contains strong positive results. On the one hand Liardet builds on the work of Klarner and Benson. He proved the following: 

\begin{thm*}[Liardet \cite{Liardet}]
	The complete growth series of a virtually abelian group is $\NG$-rational for any finite generating set.
\end{thm*}

\noindent On the other hand Grigorchuk and Nagnibeda extended results of Cannon and Gromov on rationality of numerical growth series of Gromov hyperbolic groups:

\begin{thm*}[Grigorchuk--Nagnibeda \cite{GrN}]
	The complete growth series of a Gromov hyperbolic group is $\NG$-rational for any finite generating set.
\end{thm*}

\noindent Let us also point to \cite{Changey, Coxeter_complete} for rationality of complete growth series of Coxeter groups w.r.t.\ Coxeter generating sets, and \cite{Product_complete} for stability of rationality w.r.t.\ diverse group constructions.

\bigbreak

It should be clear that a series being $R$-rational implies it being $R$-algebraic. It is also direct any morphism $R \to R'$ sends $R$-rational (resp.\ $R$-algebraic) series to $R'$-rational (resp.\ $R'$-algebraic) series. The morphisms that we consider here are the augmentation map $\epsilon\colon\NG\to \N$ and the inclusion maps  $\N\into \Z$ and $\NG\into \ZG$. Consider a group $G$ and a finite generating set $S$. For conciseness, we say a pair $(G,S)$ is $\N$-rational (resp.\ $\N G$-rational) if the associated numerical (resp.\ complete) growth series is $\N$-rational (resp.\ $\N G$-rational). Properties $\Z$- or $\Z G$-rationality, and $\N$-, $\Z$-, $\N G$ and $\Z G$-algebraicity are defined similarly. This defines a lattice of properties regarding rationality and algebraicity of both complete and numerical growth series, with multiple implications holding between them:

\begin{center}
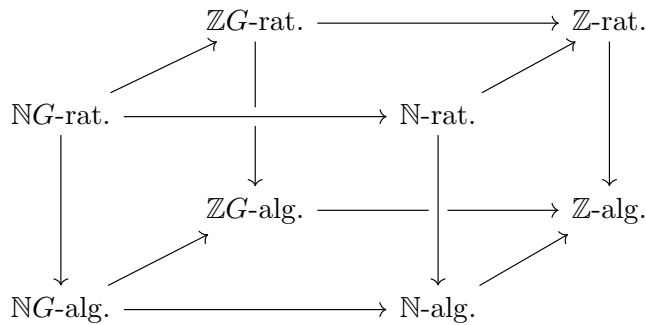

\begin{tikzcd}
    & \ZG \text{-}\rat.  \arrow{rr} \arrow{dd} & &   \Z\text{-}\rat.  \arrow{dd} \\
    \NG\text{-}\rat. \arrow[crossing over]{rr} \arrow{ur} \arrow{dd} & & \N\text{-}\rat. \arrow{ur} \\
      & \ZG\text{-}\alg.  \arrow{rr} & &  \Z\text{-}\alg. \\
    \NG\text{-}\alg. \arrow{rr} \arrow{ur} && \N\text{-}\alg. \arrow[from=uu,crossing over] \arrow{ur}
\end{tikzcd}
\captionof{figure}{Implications between different properties}
\end{center}
\noindent To the best of our knowledge, there were no counterexample that could show that the implications induced by $\epsilon$ are not equivalences. In other words, it wasn't clear if being $\N$-rational (resp. $\N$-algebraic) could imply being also $\N G$-rational (resp. $\N G$-algebraic). In the present paper, we prove this is not the case, by providing the examples of Table \ref{intro:table}.

\begin{center}
	\begin{tabular}{l|c|c|c|c}
		& \, $\N$-rational \, & \, $\N$-algebraic \, & \,$\NG$-rational\, & \,$\NG$-algebraic\,  \\
		\hline
		Virtually abelian & \cmark & \cmark & \cmark & \cmark \\
		Hyperbolic & \cmark & \cmark & \cmark & \cmark \\
		Heisenberg $H_3(\Z)$ & \cmark & \cmark & {\color{red}\xmark} & {\color{red}\xmark} \\
		$H_{5}(\Z)$ \hspace{1.8mm} (standard gen.)\! & \xmark & \xmark & \xmark & \xmark \\
		$H_{2n+1}(\Z)$ (cubical gen.)\! & \cmark & \cmark & {\color{red}\xmark} & {\color{red}\xmark} \\
		$C_2\wr \Z$ \hspace{3.4mm}(standard gen.)\! & \cmark & \cmark & {\color{red}\xmark} & \cmark \\
		$C_2\wr F_2$ \hspace{1.9mm}(standard gen.)\! & \xmark & \cmark & \xmark &  \cmark
	\end{tabular}
	\captionof{table}{Summary of our results (in red) together with some relevant known results (in black). The first three lines hold for any finite generating set.} \label{intro:table}
\end{center}

It should be noted our methods are not able to distinguish between $\Z$-rational and $\ZG$-rational groups. (Using other methods, the first author was able to do so for a natural submonoid of Thompson's group $F$.) That being said, there are no known example of groups with $\Z$-rational but non-$\N$-rational numerical growth series, or similarly with $\ZG$-rational but non-$\NG$-rational complete growth series. Therefore, we believe Table \ref{intro:table} should hold with $\N$ replaced by $\Z$.

\bigbreak

Let us summarize how those results are obtained: The obstruction to $\NG$-rationality relies on existence of dead ends of arbitrarily large depths (aka \emph{deep pockets}) - or related \say{saddle elements} - in the group $(G,S)$ of interest. Informally, a dead end of depth $D$ is an element $g\in G$ so that no longer element can be found in a ball of radius $D$ around $g$. The existence of groups with dead ends of arbitrarily large depth is not trivial, and was open for several years, until the first examples were given by \cite{Dead_ends_Lamplighter}. Here is a version of this obstruction:

\begin{thm}[Corollary \ref{sec2:crit_deep}] \label{sec0:crit}
	Let $G$ be an infinite group and $S$ a finite generating set. If $(G,S)$ has dead ends of arbitrarily large depth, then its complete growth series is not $\NG$-rational.
\end{thm}

This will apply to a large array of groups, in particular to the discrete Heisenberg group $H_3(\Z)=\la x,y\mid [x,y]\text{ central}\ra$. The growth of $H_3(\Z)$ has been under a lot of attention, culminating in results of Duchin and Shapiro \cite{Duchin_Shapiro} showing that its numerical growth series is rational (indeed $\N$-rational) w.r.t.\ any symmetric generating set. On the other hand, deep pockets were already exhibited by Warshall \cite{Warshall_Heis}. This already implies that  complete growth series for $H_3(\Z)$ is not $\NG$-rational for all generating set. We extend Warshall's results showing that \emph{all} large commutators are \emph{almost} dead ends of large depths, more precisely
\begin{thm}[Theorem \ref{almost_dead_ends_in_HZ}(a)] \label{intro:dead_ends_in_Heis}
	Fix $S$ a finite symmetric generating set for $H_3(\Z)$. There exist constants $C,M$ such that, for all $D\gg 1$ and $n\ge C\cdot D^4$, the element $g=[x,y]^n$ satisfies 
	\[ \forall h\in H_3(\Z) \text{ such that } \norm h_S \le D,\qquad \norm{gh}_S \le \norm{g}_S + M. \]
\end{thm}
\noindent  With additional arguments from Carnot-Caratheodory geometry we prove
\begin{thm}[Theorem \ref{algHeis:Heisenberg NG-alg}]
The complete growth series of the Heisenberg group $H_3(\Z)$ is not $\NG$-algebraic for all finite symmetric generating set.
\end{thm}

\bigbreak

In the case of higher Heisenberg groups $H_{2n+1}(\Z)$, Stoll \cite{Stoll} has shown that the numerical growth series of $H_5(\Z)$ w.r.t.\ the standard generating set is transcendental (i.e., not $\Z$-algebraic), while for all $n\ge 1$ the growth series of $H_{2n+1}(\Z)$ is $\N$-rational w.r.t.\  better-behaved \say{cubical} generating sets. On the other side, we establish a version of Theorem \ref{intro:dead_ends_in_Heis} for dead ends in $H_{2n+1}(\Z)$  w.r.t.\ \say{cubical-like} generating sets (Theorem \ref{almost_dead_ends_in_HZ}(b); see \S\ref{sec:higher_heis} for a definition of \say{cubical-like}). As a direct corollary the complete growth series of $H_{2n+1}(\Z)$ w.r.t.\ any cubical-like generating set is not $\NG$-rational. Work in progress of the second author shows that complete growth series of $2$-step nilpotent groups are never $\NG$-rational, unless the group is virtually abelian.

\bigbreak

Other groups of interest are wreath products $L\wr Q$. We will only be interested in generating set of the form $S_L\sqcup S_Q$. Deep pockets can be found in any group $L\wr \Z$ where $L$ itself has deep pockets (eg.\ $L=C_2$) \cite{Dead_ends_Lamplighter}. This result was recently extended to wreath products $L\wr Q$ with $Q$ abelian and any $S_Q$, and even any group and at least one generating set $S_Q$ \cite{Silva2022DeadEO}. However, as soon as the lamp group $L$ doesn't have deep pockets (eg.\ $L=\Z$), or for specific (non-abelian) base group $Q$ with specific $S_Q$ (given by $\cite{Silva2022DeadEO}$), the group $L\wr Q$ doesn't have deep pockets. This led us to introduce \say{saddle elements} (see \S2 for a definition), which then allowed us to extend our techniques to some groups without dead ends. In particular
\begin{thm}[Theorem \ref{saddle_lamplighter}]
	Let $L$ be a non-trivial group and $S_L$ a symmetric generating set. The complete growth series of the wreath product $G=L\wr \Z$ associated to the standard generating set $S_L\sqcup\{t^\pm\}$ is not $\NG$-rational.
\end{thm}
\noindent This is to be compared with
\begin{thm*}[{Bartholdi \cite[Theorem C.3]{Bartholdi}}]
	Let $L$ be a group and $S_L$ a symmetric generating set. Suppose that the complete growth series of the pair $(L,S_L)$ is $\NG$-algebraic. Then the complete growth series of the wreath products $G=L\wr \Z$ with respect to $S_L\sqcup \{t^\pm\}$ is $\NG$-algebraic.
\end{thm*}
\noindent For instance both results hold for $C_2\wr \Z$. (In this case $\NG$-algebraicity also follows from the unambiguous context-free language of geodesic representatives constructed in \cite{Language_Lamplighter}.)
\medbreak
\noindent These results extend when replacing $\Z$ by \say{tree-like} groups, i.e., groups of the form $C_2^{*p}*\Z^{*q}$ whose Cayley graphs are $m$-regular trees ($m=p+2q$). Note that, whenever $m\ge 3$, Parry \cite{Parry_Wreath} proved moreover that the numerical growth series is not rational (indeed not $\Z$-rational).

\bigbreak

\textbf{Acknowledgment.} The authors acknowledge support of the Swiss NSF grants 200020-178828 and 200020-200400. They thank Tatiana Nagnibeda, Laurent Bartholdi and Pierre de la Harpe for helpful conversations and guidances, and the anonymous referee for their comments leading to an improved presentation.
	
	\counterwithin{thm}{section}

	\section{Rationality and Algebraicity} \label{sec:NG_rat_alg}

In this section, we recall the definition of $R$-rational and $R$-algebraic formal series, where $R$ is a fixed semiring. In what follows, we will consider $R=\N G$, other relevant cases being $R=\N,\Z$ and $\Z G$. We then proceed with alternate characterizations. See \cite{handbook, Salomaa_Soittoila} for a complete treatment.

\begin{defn}[Proper systems] An \emph{algebraic system} over $R[s]$ is a system of equations
	\[ X_i = P_i(X_1,X_2,\ldots,X_n) \quad\text{for }i=1,2,\ldots,n \]
	where $P_i\in R[s]\la X_1,X_2,\ldots,X_n\ra$ are polynomials with (a priori) non-commutative coefficients and variables, i.e., finite sums of monomials $r_0X_{i_1}r_1X_{i_2}r_2\ldots X_{i_d}r_d\cdot s^k$ with $d,k\ge 0$ integers, $1\le i_j\le n$ and $r_j\in R$. A system is \emph{proper} if it contains
	\begin{itemize}
		\item no constant term (i.e., no monomial with $d=k=0$);
		\item no monomial with $d=1$ and $k=0$.
	\end{itemize}
	A system is \emph{linear} if furthermore all monomials satisfy $d=0$, or $d=1$ and $r_1=1_R$.
\end{defn}
\begin{defn}
	A formal series $S=\sum_{n=0}^\infty\,a_n\cdot s^n\in R[[s]]$ is \emph{proper} (or \emph{quasiregular}) if its constant coefficient is zero, that is $a_0=0$.
\end{defn}
\bigbreak
The main motivation for those two definitions is the following result:
\begin{thm}[{See \cite[Theorem IV.1.1]{Salomaa_Soittoila}}]
	Every proper algebraic system admits a unique solution $(S_1,S_2,\ldots,S_n)\in R[[s]]^n$ consisting of proper formal series. 
\end{thm}
Finally, we say a proper series is $R$-rational (resp. $R$-algebraic) if it is solution to a proper linear (resp. algebraic system). More generally
\begin{defn}[$R$-rational and $R$-algebraic series]
A formal series $ S \in R[[s]]$ is \emph{$R$-rational} (resp.\ \emph{$R$-algebraic}) if $S-a_0$ is the first component of a proper solution $(S_1,S_2,\ldots,S_n)\in R[[s]]$ to a proper linear (resp.\ algebraic) system over $R[s]$.
\end{defn}
\bigbreak

Let us give another less technical definition of rational series over non-commutative semirings. We first need the notion of quasi-inverse:

\begin{defn}
Given a \emph{proper} series $S(s)\in R[[s]]$, we define its \emph{quasi-inverse} as
\[ S(s)^* = \sum_{n\ge 0} S(s)^n,\]
where $S(s)^n$ is the multiplication of $S(s)$ with itself $n$ times and $S(s)^0=1$.
\end{defn}

The name quasi-inverse comes from the \say{equality} $(1-S)S^*=1$. This (linear) equation shows that, if $S$ itself is rational, then $S^*$ is rational too. Reciprocally, all rational series can be obtained from polynomials using addition, multiplication and quasi-inversion:

\begin{thm}[{See \cite[Theorems II.1.2 and II.1.4]{Salomaa_Soittoila}}]\label{sec1:closed_quasiinv}
	The set of $R$-rational series is the smallest semiring containing $R[s]$ and closed under quasi-inversion (for proper series).
\end{thm}

Rational and algebraic series can also be characterized in terms of automata and languages.

\begin{defn} An \emph{$R[s]$-automaton} is a finite oriented graph $A=(V,E)$ with an edge labeling $p\colon E\to R[s]$. A series $S(s)\in R[[s]]$  is \emph{recognized by the automaton} $A$ if there exist two initial and terminal vertices $I,T\in V$ such that
\[ S(s) = \sum_{\gamma \in P(I,T)} p(\gamma),\]
where $P(I,T)$ is the set of all oriented paths from $I$ to $T$, and $p(e_1\ldots e_\ell) := p(e_1)\ldots p(e_\ell)$.
\end{defn}
\begin{thm}[{Kleene-Schützenberger \cite{Schutzenberger_Fund}, see also \cite[Theorem 3.2.5]{handbook}}] \label{sec1:N_automaton}
A formal series is $R$-rational if and only if it is recognized by an $R[s]$-automaton.
\end{thm}
\begin{rmk}\label{sec1:rem}  ~
\begin{enumerate}[leftmargin=8mm, label=(\alph*)]
	\item In case $R=\N G$, we may \emph{and will} assume that all labels are of the form $p(e)=w(e)\cdot s^{l(e)}$ for some $w(e)\in G$ and $l(e)\in\N$, at the cost of allowing multiple edges between vertices.
	\item If furthermore $S(s)$ is the complete growth series, then $\|w(\gamma)\|=l(\gamma)$ holds for all path $\gamma \in P(I,T)$ (with the obvious notations)
\end{enumerate}
\end{rmk}
\begin{rmk} \label{sec1:rmk_compu}
	The computational advantages of $\NG$-rational series over $\ZG$-rational series are made clearer through the formalism of $R[s]$-automata. Given an $\NG[s]$-automaton for a complete growth series, we can pick a path at random, compute the associated term $g\cdot s^l$, and get a certificate that the element $g$ has length $l$. In contrast, if the automaton had labels in $\ZG[s]$, the term $g\cdot s^l$ might cancel out with other terms, so that we need to compute all terms with the same exponent $s^l$ before having a certificate $g$ has indeed length $l$.
\end{rmk}

Theorem \ref{sec1:N_automaton} admits a generalization linking $R$-algebraic series and weighted pushdown automata (see \cite[Section 7.5.3]{handbook}). For complete growth series, this can be reformulated as
\begin{thm} \label{sec5:CF_cross_section}
	Let $G$ be a group and $S$ a generating set. The associated complete growth series is $\N G$-algebraic if and only if there exists an unambiguous context-free language $\Lc\subset (G\times\N)^*$ such that the evaluation map $\ev\colon \Lc\to G\times\N$ is injective and its image is given by
	\[ \ev(\Lc) = \{ (g,\norm g_S) : g\in G\}. \]
\end{thm}
	\section{Dead ends and saddle elements}

Let us recall the notion of dead ends introduced in \cite{Bogopolskii}, and introduce some related notions:
\begin{defn}
	Consider $(G,\norm{\,\cdot\,})$ a group endowed with a left-invariant metric (e.g., a finitely generated group endowed with a word metric, or a nilpotent Lie group endowed with a Carnot-Cartheodory metric) and let $D>0$ and $M\ge 0$ be real numbers. An element $g\in G$ is
	\begin{itemize}[leftmargin=8mm]
		\item a \emph{dead end of depth} at least $D$ if 
		\[ \forall h\in G \text{ such that }\norm h\le D, \quad \norm{gh}\le \norm g.\]
		\item an \emph{almost dead end of depth} at least $D$ and \emph{margin $M$} if
		\[ \forall h\in G \text{ such that }\norm h\le D, \quad \norm{gh}\le \norm g + M.\]
		\item a \emph{saddle element of depth} at least $D$ if
		\[ \forall h\in G \text{ such that }\norm h\le D, \quad \norm{gh}=\|gh^{-1}\|.\]
	\end{itemize}
	We say that $(G,\norm{\,\cdot\,})$ has \emph{deep pockets} if it possesses dead ends of arbitrarily large depths.
\end{defn}
\medbreak
\begin{rmk}
	Finite (and more generally compact) groups have deep pockets: element of maximal length are dead ends of infinite depth. The notion of deep pockets gets way more interesting when looking at non-bounded groups, as all elements have finite depth.
\end{rmk}
\medbreak
To the best of our knowledge, the notion of almost dead ends is new, even though almost dead ends appear implicitly in Warshall's work \cite{Warshall_Latt, Warshall_Heis}. The main interest of almost dead ends is that the existence of almost dead ends of large depth is preserved by rough isometries.
\begin{defn}
	Let $X,Y$ be two metric spaces. A function $f\colon X\to Y$ is a \emph{rough isometry} {\normalfont(}or $(1,K)$-quasi-isometry{\normalfont)} if there exists some constant $K\ge 0$ such that
	\begin{itemize}[leftmargin=8mm]
		\item for all $x,x'\in X$, we have
		\[ d_X(x,x') - K \le d_Y\big(f(x),f(x')\big) \le d_X(x,x') + K.  \]
		\item for all $y\in Y$, there exists $x\in X$ such that $d_Y\big(y,f(x)\big)\le K$. \vspace*{2mm}
	\end{itemize}
\end{defn}
\noindent\textbf{Remark.} Even though this might not be obvious from the definition, this defines an equivalence relation on metric spaces, in particular it is symmetric.

\begin{lem}[{Compare with \cite[Proposition 7]{Warshall_Latt}}] \label{sec1:rough_dead_ends}
	Let $(G,d_G)$ and $(H,d_H)$ be two roughly isometric metric group. Suppose that there exists some $M\ge 0$ such that $G$ contains almost dead ends of arbitrarily large depth and fixed margin $M$, then the same is true in $H$.
\end{lem}
\begin{proof}
	Let $f\colon G\to H$ be a rough isometry. Up to translation we may suppose that $f(e_G)=e_H$. Let $g\in G$ be an almost dead end of margin $M$ and depth $D$. We prove that $f(g)$ is an almost dead end. For any $h\in H$ s.t.\ $\| h\| \le D-2K$, there exists $x\in G$ s.t.\ $d\big(f(x), f(g) h\big)\le K$ hence
	\[ d(g,x) \le d\big(f(g),f(x)\big)+K \le d\big(f(g),f(g) h\big)+d\big(f(g) h,f(x)\big) + K \le D.\]
	By hypothesis $g$ is an almost dead end of depth $D$ so that $\norm x\le \norm g+M$ hence
	\[ \norm{f(g) h} \le \norm{f(x)} + d\big(f(x), f(g) h\big) \le \big(\!\norm x+K\big) + K \le \norm g + M + 2K \le \norm{f(g)} + M + 3K. \]
	This means that $f(g)$ is an almost dead end of margin $M+3K$ (which is fixed) and depth $D-2K$ (which can be made arbitrarily large for well chosen $g\in G$).
\end{proof}
\medbreak
Usually we should think of the situation where the depth $D$ is way larger than the margin $M$. In case of integer-valued metric and $D\gg M$, the following key lemma ensures that we can \say{promote} almost dead ends into genuine dead ends of comparatively large depth.
\begin{lem}[{Warshall \cite[Proposition 6]{Warshall_Latt}}] Let $X$ be a metric space, $f\colon X\to \N$ a function, and $D,M\in\N$. Suppose that there exists $x\in X$ such that
	\[ f(x') < f(x)+M \text{ for all }x'\in B_D(x).\]
	Then there exists $x'\in B_D(x)$ such that $f$ reaches a maximum on $B_{\frac DM}(x')$ at $x'$.
\end{lem}

\bigbreak
Let us now relate dead ends with the rationality of complete growth series
\begin{thm}[Pumping lemma + $\varepsilon$]\label{sec2:crit_gene}
	Let $G$ be a group  and $S$ a finite generating set. Suppose that the complete growth series of $(G,S)$ is $\NG$-rational. Then there exists $D>0$ such that all but finitely many $g\in G$ can be written as $g=uv$ with $v\ne e$, $\norm v \le D$ and
		\[ \norm{uv^n} = \norm{u} + n \tau \qquad\text{for all }n\ge 0 \]
	for a certain $\tau\in\N\setminus\{0\}$. Moreover 
	\[ \tau=\tau(v) = \lim_{n\to\infty} \frac{\norm{v^n}}n \]
	is the \emph{translation number} of $v$.
\end{thm}
\begin{proof}
	As the complete growth series is $\NG$-rational, it is recognized by an automaton $(V,E)$. We take $D=2\abs{V}\max_{e\in E}\norm{w(e)}$ (with $p(e)=w(e)s^{l(e)}$ as in Remark \ref{sec1:rem}). Note that each element $g\in G$ corresponds to a path $\gamma_g\in P(I,T)$. Moreover, all but finitely many $g$ correspond to paths containing cycles. In that case, we decompose $\gamma_g$ as a concatenation $\alpha\beta\delta$ with $\beta$ a \emph{non-empty} simple cycle and $\delta$ a simple path. Among the $\abs V+1$ last visited vertices, at least one is repeated. The loop between the two last visits of this vertex is $\beta$, and everything coming afterward is $\delta$. In particular the path $\beta\delta$ is formed of at most $\abs V$ edges.
\medbreak
Take $u=w(\alpha\delta)$, $v=w(\delta)^{-1}w(\beta)w(\delta)$ and $\tau=l(\beta)\in \N$. Note that $\norm v \le D$ as $v$ is a product of at most $2\abs V$ edge labels $w(e)$. Moreover $\alpha\beta^n\delta\in P(I,T)$ for all $n$, so that we get terms $w(\alpha\beta^n\delta)\cdot s^{l(\alpha\beta^n\delta)}$ in the complete growth series. It follows that
\[ \norm{uv^n} = \norm{w(\alpha\beta^n\delta)} = l(\alpha\beta^n\delta) = l(\alpha\delta)+n\tau = \norm{w(\alpha\delta)} +n\tau = \norm{u}+n\tau. \]
\indent Finally, let us look at $\tau$. The previous equality implies that $\tau\ne 0$, otherwise we would get infinitely many elements inside the sphere of radius $\norm{u}$. 
Moreover the triangle inequality gives
\[ \abs{\tau-\frac{\norm{v^n}}n} = \abs{\frac{\norm{uv^n}-\norm{u}-\norm{v^n}}n} \le \frac{2\norm u}n \xrightarrow{n\rightarrow \infty} 0. \]
so that $\tau=\lim_{n\to\infty}\frac{\norm{v^n}}n$ is indeed the translation number of $v$.
\end{proof}

Theorem \ref{sec0:crit} announced in the introduction follows as a quick corollary:
\begin{co}\label{sec2:crit_deep}
	Let $G$ be a group and $S$ a finite generating set. If either
	\begin{enumerate}[leftmargin=8mm, label={\normalfont(\roman*)} ]
		\item $G$ is infinite and $(G,S)$ has deep pockets, or
		\item $(G,S)$ contains an infinite sequence of distinct saddle elements of increasing depth,
	\end{enumerate}
then its complete growth series of $(G,S)$ is not $\N G$-rational.
\end{co}
\begin{proof}
We argue by contraposition. Suppose that $(G,S)$ is infinite and its complete growth series is $\N G$-rational. By Theorem \ref{sec2:crit_gene} all but finitely many elements can be written as $g=uv$ with $\norm v\le D$ and $\norm{uv^n} = \norm{u} + n \tau$  for all $n\ge 0$. It follows that
\[ \norm{gv}=\|uv^2\|=\norm{g}+\tau > \norm g > \norm g-\tau = \norm{u} = \|gv^{-1}\|. \]
\begin{enumerate}[leftmargin=8mm, label=(\roman*) ]
	\item In particular $g$ cannot be a dead end of depth $\ge D$. As $G$ is infinite, all elements have finite depth, in particular the finitely many elements left out by this argument have finite depth. It follows that $(G,S)$ doesn't have deep pockets.
	
	\item In particular $g$ cannot be a saddle element of depth $\ge D$. As only finitely many elements are left out, $G$ cannot contain an infinite sequence with the required condition. \qedhere
\end{enumerate}
\end{proof}
	\section{Dead ends in Heisenberg groups} \label{sec:Heis}

We reprove one of Warshall's results, in order to retrieve the finer conclusions of Theorem \ref{intro:dead_ends_in_Heis}, which will then allow us to find dead ends in higher Heisenberg groups. 

\subsection{Models for $H_3(\R)$ and $H_3(\Z)$}

Elements of $H_3(\R)$ can be seen as equivalence classes of absolutely continuous curves in $\R^2$ starting at $\mathbf 0$, with concatenation as operation (see Figure \ref{fig:conc}). We associate to each curve $\gamma$
\begin{itemize}[leftmargin=8mm]
	\item Its second endpoint $\hat\gamma\in\R^2$. Note that $\gamma\mapsto \hat\gamma$ maps $H_3(\R)$ to its abelianization $\R^2$.
	\item Its \emph{balayage area}\footnote{\say{Balayer} is french for \say{to sweep}. $A(\gamma)$ is the area swept by the moving segment from $\mathbf 0$ to $\gamma(t)$.} $A(\gamma) =  \int (\gamma_x\gamma_y'-\gamma_x'\gamma_y)\cdot \dif t \in\R$.
\end{itemize}
Two curves $\gamma_1,\gamma_2$ represent the same element of $H_3(\R)$ if $\hat\gamma_1=\hat\gamma_2$ and $A(\gamma_1)=A(\gamma_2)$. 


\begin{prop}
	Given two paths $g,h$, their concatenation $gh$ has parameters
	\begin{align*}
	 \widehat{gh} & = \hat g+\hat h, \\
	 A(gh) & = A(g) + A(h) + \frac12 \det(\hat g,\hat h).
	\end{align*}\vspace*{-5mm}

	\noindent In particular, the operation \say{concatenation} is well-defined on $H_3(\R)$. With the constant curve as neutral element, and curves traveled backward as inverses, this defines a group.
\end{prop}
\begin{center}
	\begin{tikzpicture}[scale=1.28]
			\clip (-.8,-.5) rectangle (2.55,3.15);
			\fill[Purple!40, opacity=.5] (0,0)
				to[out=-10, in=-110, looseness=1.5] (1.2,.6)
				to[out=70, in=160, looseness=2] (2.2,.8);
			\draw[Purple, thick, -latex] (0,0)
				to[out=-10, in=-110, looseness=1.5] (1.2,.6)
				to[out=70, in=160, looseness=2] (2.2,.8);
			\draw[dashed, Purple, thick, latex-] (0,0) -- (2.2,.8);
			
			\node[Purple] at (.7,.11) {\scriptsize$+1$};
			\node[Purple] at (1.55,.8) {\scriptsize$-1$};
		
			\fill[orange!40, opacity=.5] (0,0)
				to [out=70, in=-60, looseness=1.5] (0,1.3)
				to [out=120, in=65, looseness=1.5] (-.4,1.2)
				to [out=-115, in=-90, looseness=1.2] (0,1.3)
				to [out=90, in=65, looseness=2] (-.5,1.3)
				to [out=-115, in=-150, looseness=2] (.2,.9)
				to [out=30, in=-90, looseness=1.5] (1,2);
			\fill[orange!40]	(0,1.3)
				to [out=120, in=65, looseness=1.5] (-.4,1.2)
				to [out=-115, in=-90, looseness=1.2] (0,1.3);
			\draw[orange, thick, -latex] (0,0)
				to [out=70, in=-60, looseness=1.5] (0,1.3)
				to [out=120, in=65, looseness=1.5] (-.4,1.2)
				to [out=-115, in=-90, looseness=1.2] (0,1.3)
				to [out=90, in=65, looseness=2] (-.5,1.3)
				to [out=-115, in=-150, looseness=2] (.2,.9)
				to [out=30, in=-90, looseness=1.5] (1,2);
			\draw[dashed, orange, thick, latex-] (0,0) -- (1,2);
			
			\node[Orange] at (-.2,.96) {\scriptsize$+1$};
			\node[Orange] at (-.2,1.25) {\scriptsize$+2$};
			
			\node[circle, fill=black, inner sep=1.5pt, label=below:{\footnotesize$\mathbf 0$}] at (0,0) {};
			\node[circle, fill=Purple, inner sep=1pt, label=below:{\color{Purple}$\hat g$}] at (2.2,.8) {};
			\node[circle, fill=orange, inner sep=1pt, label=above:{\color{orange}$\hat h$}] at (1,2) {};
		\end{tikzpicture} \hspace*{2mm}
		\begin{tikzpicture}[scale=1.28]
			\clip (-.35,-.5) rectangle (3.55,3.15);
			
			\fill[TealBlue!40, opacity=.5] (0,0) -- (2.2,.8) -- (3.2,2.8) -- cycle;
			\draw[dashed, TealBlue, thick, latex-] (0,0) -- (3.2,2.8);
			
			\fill[Purple!40, opacity=.5] (0,0)
				to[out=-10, in=-110, looseness=1.5] (1.2,.6)
				to[out=70, in=160, looseness=2] (2.2,.8);
			\draw[Purple, thick, -latex] (0,0)
				to[out=-10, in=-110, looseness=1.5] (1.2,.6)
				to[out=70, in=160, looseness=2] (2.2,.8);
			\draw[dashed, Purple, thick] (0,0) -- (2.2,.8);
			\node[Purple] at (.7,.11) {\scriptsize$+1$};
			\node[Purple] at (1.55,.8) {\scriptsize$-1$};
			
			\begin{scope}[shift={(2.2,.8)}]
				\fill[orange!40, opacity=.5] (0,0)
					to [out=70, in=-60, looseness=1.5] (0,1.3)
					to [out=120, in=65, looseness=1.5] (-.4,1.2)
					to [out=-115, in=-90, looseness=1.2] (0,1.3)
					to [out=90, in=65, looseness=2] (-.5,1.3)
					to [out=-115, in=-150, looseness=2] (.2,.9)
					to [out=30, in=-90, looseness=1.5] (1,2);
				\fill[orange!40]	(0,1.3)
					to [out=120, in=65, looseness=1.5] (-.4,1.2)
					to [out=-115, in=-90, looseness=1.2] (0,1.3);
				\draw[orange, thick, -latex] (0,0)
					to [out=70, in=-60, looseness=1.5] (0,1.3)
					to [out=120, in=65, looseness=1.5] (-.4,1.2)
					to [out=-115, in=-90, looseness=1.2] (0,1.3)
					to [out=90, in=65, looseness=2] (-.5,1.3)
					to [out=-115, in=-150, looseness=2] (.2,.9)
					to [out=30, in=-90, looseness=1.5] (1,2);
				\draw[dashed, orange, thick] (0,0) -- (1,2);
				\node[Orange] at (-.2,.96) {\scriptsize$+1$};
				\node[Orange] at (-.2,1.25) {\scriptsize$+2$};
			\end{scope}
			
			\node[circle, fill=black, inner sep=1.5pt, label=below:{\footnotesize$\mathbf 0$}] at (0,0) {};
			\node[circle, fill=Purple, inner sep=1pt, label=below:{\color{Purple}$\hat g$}] at (2.2,.8) {};
			\node[circle, fill=orange, inner sep=1pt] at (3.2,2.8) {};
			\node at (3.2,3) {\footnotesize$\hat g+\hat h$};
		
			\node[TealBlue] at (2,1.25) {\footnotesize$+1$};
		\end{tikzpicture}
		\begin{tikzpicture}[scale=1.28]
		\clip (-.35,-.5) rectangle (3.55,3.15);
		
		\begin{scope}[shift={(2.2,.8)}]
		\fill[TealBlue!40, opacity=.5] (-2.2,-.8)
			to[out=-10, in=-110, looseness=1.5] (-1,-.2)
			to[out=70, in=160, looseness=2] (0,0)
			to [out=70, in=-60, looseness=1.5] (0,1.3)
			to [out=120, in=65, looseness=1.5] (-.4,1.2)
			to [out=-115, in=-90, looseness=1.2] (0,1.3)
			to [out=90, in=65, looseness=2] (-.5,1.3)
			to [out=-115, in=-150, looseness=2] (.2,.9)
			to [out=30, in=-90, looseness=1.5] (1,2);
		\fill[TealBlue!40]	(0,1.3)
			to [out=120, in=65, looseness=1.5] (-.4,1.2)
			to [out=-115, in=-90, looseness=1.2] (0,1.3);
		\begin{scope}
			\clip (-2.2,-.8) -- (1,2) -- (.2,.9) -- cycle;
			\fill[TealBlue!40] (0,0)
				to [out=70, in=-60, looseness=1.5] (0,1.3)
				to [out=120, in=65, looseness=1.5] (-.4,1.2)
				to [out=-115, in=-90, looseness=1.2] (0,1.3)
				to [out=90, in=65, looseness=2] (-.5,1.3)
				to [out=-115, in=-150, looseness=2] (.2,.9);
		\end{scope}
		\draw[TealBlue, thick, -latex] (-2.2,-.8)
			to[out=-10, in=-110, looseness=1.5] (-1,-.2)
			to[out=70, in=160, looseness=2] (0,0)
			to [out=70, in=-60, looseness=1.5] (0,1.3)
			to [out=120, in=65, looseness=1.5] (-.4,1.2)
			to [out=-115, in=-90, looseness=1.2] (0,1.3)
			to [out=90, in=65, looseness=2] (-.5,1.3)
			to [out=-115, in=-150, looseness=2] (.2,.9)
			to [out=30, in=-90, looseness=1.5] (1,2);
			\node[TealBlue] at (-.2,1.25) {\scriptsize$+2$};
			\node[TealBlue] at (-.02,.94) {\scriptsize$+2$};
			\node[TealBlue] at (-.41,.98) {\scriptsize$+1$};
		\end{scope}
		\node[TealBlue] at (.7,.11) {\scriptsize$+1$};
		\node[TealBlue] at (1.55,.8) {\scriptsize$0$};
		
		\draw[dashed, TealBlue, thick, latex-] (0,0) -- (3.2,2.8);
		\node[circle, fill=black, inner sep=1.5pt, label=below:{\footnotesize$\mathbf 0$}] at (0,0) {};
		\node[circle, fill=TealBlue, inner sep=1pt] at (3.2,2.8) {};
		\node at (3.2,3) {\footnotesize$\hat g+\hat h$};
		
		\node[TealBlue] at (2,1.25) {\footnotesize$+1$};
		\end{tikzpicture}
	
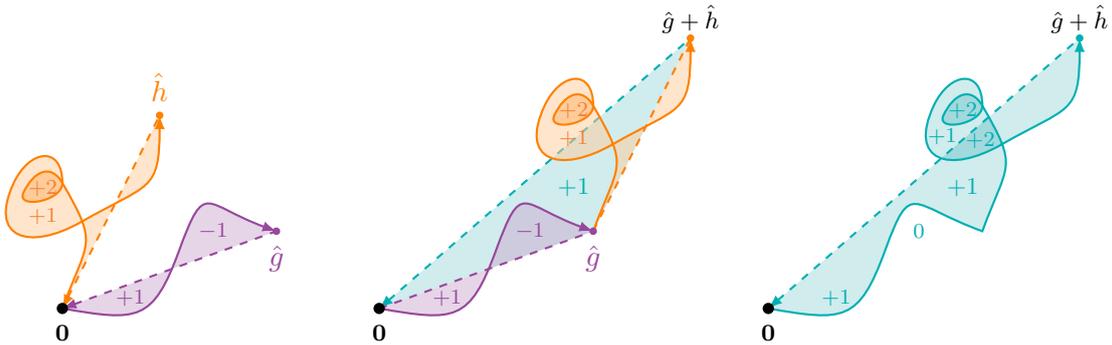
\captionof{figure}{Concatenation of two elements $g,h$ and a few winding numbers.} \label{fig:conc}
\end{center}
\noindent An alternate formula for $A(\gamma)$ helps to better understand this result. We define a closed curve $\gamma_c$ by concatenating $\gamma$ with the  segment from $\hat\gamma$ to $\mathbf 0$, and a function $w_\gamma\colon\R^2\setminus\mathrm{Im}(\gamma_c)\to\Z$ defined as $w_\gamma(x,y)=$ winding number of $\gamma_c$ around $(x,y)$ (see Figure \ref{fig:conc}). Then
	\[ A(\gamma) = \iint_{\R^2} w_\gamma(x,y) \cdot \dif x \,\dif y \in\R. \vspace*{1mm} \]	

\newpage

\noindent The discrete Heisenberg group $H_3(\Z)$ can be seen as the subgroup generated by the unit segments $x,y$ from the origin $\mathbf 0$ to $(1,0), (0,1)$ respectively. Moreover, we denote $z=[x,y]$.

\subsection{$\CC$-metrics on $H_3(\R)$}

We recall the definition of the Carnot-Caratheodory metric on $H_3(\R)$. First
\begin{defn}[Minkowski norm]
	Let $L\subset \R^n$ be a centrally symmetric (full dimensional) convex polytope. We define the \emph{Minkowski norm} on $\R^n$ as $\norm v_L = \min\{ \lambda\ge 0 : v\in \lambda L\}$.
\end{defn}
\begin{defn}[$\CC$-metric on $H_3(\R)$] \label{app:def:CC_metric}
	Given $L\subset \R^2$ centrally symmetric convex polygon, we define the Carnot-Caratheodory metric -- or $\CC$-metric for short -- on $H_3(\R)$ as
	\[ \norm g_\CC = \min\left\{ \norm \gamma_L \;\big|\; \gamma\text{ is a path representing }g\right\}, \]
	where $\norm{\gamma}_L \coloneqq \int \norm{\gamma'(t)}_L \cdot \dif t$. Any path reaching this minimum is a $\CC$-geodesic.
\end{defn}

This norm is particularly nice among subFinsler metrics as it is homogeneous:
\begin{prop}[Dilation] \label{prop:dilation}
	Given $\lambda\in\mathbb R$, we define the \emph{dilation} $\delta_\lambda\colon H_3(\R)\to H_3(\R)$ as the map sending each path to the image of this path by an homothety of ratio $\lambda$ {\normalfont(}centered at $\mathbf 0${\normalfont)}.  The map $\delta_\lambda$ is an automorphism {\normalfont(}for $\lambda\ne 0${\normalfont)} satisfying
	\[ \widehat{\delta_\lambda g} = \lambda \hat g,\quad A(\delta_\lambda g)=\lambda^2A(g) \quad\text{and}\quad \norm{\delta_\lambda g}_\CC = \abs{\lambda}\cdot \norm g_\CC. \]
\end{prop}

\bigbreak

$\CC$-metrics on $H_3(\R)$ can be described precisely. First a classical result of Busemann on the isoperimetric problem in the plane with Minkowski metric:
\begin{thm}[Busemann \cite{Busemann}]
	Fix $L\subset\R^2$ centrally symmetric convex polygon. Among closed curves, the maximal ratio of enclosed area to the square of the perimeter is achieved (uniquely up to scaling) for the isoperimetrix $I$, i.e., the boundary of the rotation by $\pm\frac\pi2$ of the polar dual
	\[ L^* = \left\{ x\in \R^2 \;\big|\; \forall y\in L,\; \la x;y\ra \le 1 \right\}. \vspace*{-2mm}\]
	In other words, for any closed curve $\gamma$, we have $\frac{\norm\gamma_L^2}{A(\gamma)} \ge \frac{\norm I_L^2}{A(I)}$. 
\end{thm}
\begin{center}
	\begin{tikzpicture}[scale=.37]
	\clip (-6.2,-5) rectangle (6.2,5);
	\draw[thick, -latex] (-6,0) -- (6,0);
	\draw[thick, -latex] (0,-5) -- (0,5);
	\draw[fill=black] (2.97,-.2) rectangle (3.03,.2);
	\draw[fill=black] (-.2,2.97) rectangle (.2,3.03);
	\draw[thick, Turquoise, fill=Turquoise, fill opacity=.4] (4,0) -- (6,4) -- (0,4) -- (-4,0) -- (-6,-4) -- (0,-4) -- (4,0);
\end{tikzpicture}
\hspace*{5mm}
\begin{tikzpicture}[scale=.37]
	\clip (-6.2,-5) rectangle (6.2,5);
	
	\draw[thick, dashed, magenta] (0,0) -- (6,4);
	\draw[thick, magenta] (3,-9/4) -- (-1,15/4);
	\draw[Turquoise, fill=Turquoise] (6,4) circle (4pt);
	
	\draw[thick, dashed, Purple] (0,0) -- (0,4);
	\draw[thick, Purple] (-5,9/4) -- (5,9/4);
	\draw[Turquoise, fill=Turquoise] (0,4) circle (4pt);
	
	\draw[thick, dashed, pink] (0,0) -- (4,0);
	\draw[thick, pink] (9/4,-5) -- (9/4,5);
	\draw[Turquoise, fill=Turquoise] (4,0) circle (4pt);
	
	\begin{scope}[rotate=180]
		\draw[thick, dashed, magenta] (0,0) -- (6,4);
		\draw[thick, magenta] (3,-9/4) -- (-1,15/4);
		\draw[Turquoise, fill=Turquoise] (6,4) circle (4pt);
		
		\draw[thick, dashed, Purple] (0,0) -- (0,4);
		\draw[thick, Purple] (-5,9/4) -- (5,9/4);
		\draw[Turquoise, fill=Turquoise] (0,4) circle (4pt);
		
		\draw[thick, dashed, pink] (0,0) -- (4,0);
		\draw[thick, pink] (9/4,-5) -- (9/4,5);
		\draw[Turquoise, fill=Turquoise] (4,0) circle (4pt);
	\end{scope}
	
	\fill[magenta, opacity=.2] (9/4,-9/4) -- (9/4,-9/8) -- (0,9/4) -- (-9/4,9/4) -- (-9/4,9/8) -- (0,-9/4) -- cycle;
	
\end{tikzpicture}
\hspace*{2mm}
\begin{tikzpicture}[scale=.37]
	\clip (-4,-5) rectangle (4,5);
	
	\draw[thick, LimeGreen, rotate=-90, -latex] (9/4,-9/4) -- (9/4,-9/8) -- (0,9/4) -- (-9/4,9/4) -- (-9/4,9/8) -- (0,-9/4) -- (9/4,-9/4);
\end{tikzpicture}
	
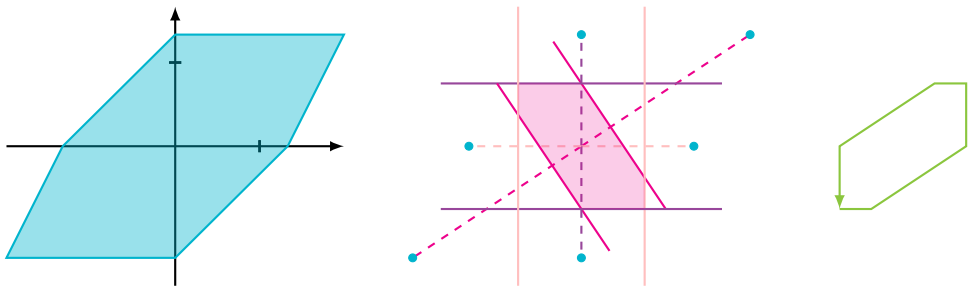
\captionof{figure}{Examples of $\color{Turquoise}L$, $\color{magenta}L^*$ and $\color{LimeGreen}I$.}
\end{center}
\noindent Note that, whenever $L$ is polygonal, then $I$ is polygonal with sides parallels to vertices of $L$.

\bigbreak

This was generalized for all $\CC$-geodesics in a beautiful paper of Duchin and Mooney:
\begin{thm}[{Duchin-Mooney \cite[Structure Theorem]{Mooney_Duchin}}] \label{thm:stru} $\CC$-geodesics split into two classes:
	\begin{itemize}[leftmargin=8mm]
		\item \emph{Regular geodesics} which follows a portion of a dilate of the isoperimetrix.
		\item \emph{Unstable geodesics} for which all tangent directions $\gamma'(t)$ lie in a common positive cone spanned by two consecutive vertices of $L$. These coincides with geodesics in $(\R^2,\norm{\,\cdot\,}_L)$. 
	\end{itemize}
Moreover all such paths are $\CC$-geodesics.
\end{thm}

\newpage

\begin{center}
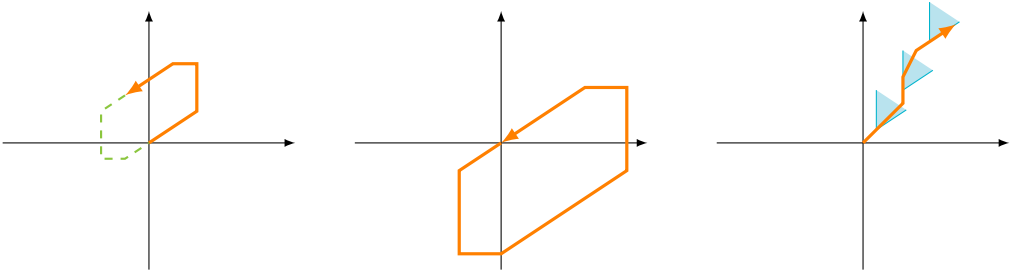

	\begin{tikzpicture}[scale=.35]
	\draw[-latex] (-5.5,0) -- (5.5,0);
	\draw[-latex] (0,-4.8) -- (0,5);
	\begin{scope}[scale=.2, shift={(0,6)}]
		\draw[very thick, orange, -latex] (0,-6) -- (9,0) -- (9,9) -- (9/2,9) -- (-9/2,3);
		\draw[thick, LimeGreen, dashed] (-9/2,3) -- (-9,0) -- (-9,-9) -- (-9/2,-9) -- (0,-6);
	\end{scope}
\end{tikzpicture}\qquad
\begin{tikzpicture}[scale=.35]
	\draw[-latex] (-5.5,0) -- (5.5,0);
	\draw[-latex] (0,-4.8) -- (0,5);
	\begin{scope}[scale=.35, shift={(9/2,-3)}]
		\draw[very thick, orange, -latex] (-9/2,3) -- (-9,0) -- (-9,-9) -- (-9/2,-9) -- (9,0) -- (9,9) -- (9/2,9) -- (-9/2,3);
	\end{scope}
\end{tikzpicture} \qquad
\begin{tikzpicture}[scale=.35]
	\draw[-latex] (-5.5,0) -- (5.5,0);
	\draw[-latex] (0,-4.8) -- (0,5); 
	
	\draw[Turquoise, thin, scale=.5, shift={(1,1)}, fill=Turquoise!30] (0,3) -- (0,0) -- (2.25,1.5);
	\draw[Turquoise, thin, scale=.5, shift={(3,4)}, fill=Turquoise!30] (0,3) -- (0,0) -- (2.25,1.5);
	\draw[Turquoise, thin, scale=.5, shift={(5,7.67)}, fill=Turquoise!30] (0,3) -- (0,0) -- (2.25,1.5);
	
	\draw[very thick, orange, scale=.5, -latex] (0,0) -- (3,3) -- (3,5) -- (4,7) -- (7,9);
\end{tikzpicture}
	
	
	\captionof{figure}{A few $\CC$-geodesics of both types for the previous choice of $L$.}
\end{center}
\begin{rmk}
	It follows that, for any $g\in [H,H]$, all $\CC$-geodesics follow the same re-scaled copy of $I$. Specifically, if $g=z^A$ (that is, $A(g)=A$), the scaling factor should be $\sqrt{\abs{A}/\abs{A(I)}}$. Only the starting point along the curve is up to choice. In particular,
	\[ \norm{g}_\CC = \sqrt{\abs A}\norm z_\CC = \sqrt{\frac{\abs{A(g)}}{\abs{A(I)}}}\norm I_L. \]
\end{rmk}

\subsection{Almost dead ends in $H_3(\R)$ for $\CC$-metrics}

It should first be noted that $H=H_3(\R)$ doesn't have any dead end\footnote{There is another definition of dead ends in pointed geodesic metric spaces. A dead end is an element $g$ such that no geodesic leading to $g$ can be extended into a geodesic for a longer element. With this definition, non-trivial commutators in $H_3(\R)$ are indeed dead ends.}. Indeed, it can easily be seen from the Duchin-Mooney Structure Theorem that $\CC$-geodesics for elements $g\notin [H,H]$ can be extended into geodesics for longer elements. Moreover, if $g=z^A$ with $A\ne 0$. For any $D>0$, we can take $h=z^{\pm \varepsilon}$ which has $\CC$-length $\sqrt{\epsilon}\norm{z}_\CC\le D$. We can chose the sign such that the equality ${\abs{A\pm \varepsilon}} > {\abs A}$ is true. 
It follows that
\[ \norm{gh}_\CC= \norm{z^{A\pm \varepsilon}}_\CC = \sqrt{\abs{A\pm \varepsilon}} \norm{z}_\CC > \sqrt{\abs A}\norm{z}_\CC =\norm g_\CC \]
for the right choice of sign. This means we can only hope for almost dead ends in $(H_3(\R),\norm{\,\cdot\,}_\CC)$, which happens to be sufficient to draw conclusions in the word metric case.
\begin{prop} \label{almost_dead_ends_in_HR}
	Fix $L\subset \R^2$ centrally symmetric convex polygon. There exist constants $C,M$ depending on $L$ such that, for any $D\gg 1$ and $n\ge C\cdot D^4$, the element $g=z^n=[x,y]^n$ satisfies
	\[ \forall h\in H_3(\R)\text{ s.t. } \norm h_\CC \le D, \quad \norm{gh}_\CC \le \norm g_\CC - \|\hat h\|_L +M. \]
\end{prop}
\begin{proof}
Fix $g,h\in H_3(\R)$ as in the statement. There exist consecutive vertices $\hat s_i,\hat s_{i+1}$ of $L$ such that $\hat h$ lives in the positive cone spanned by $\hat s_i,\hat s_{i+1}$, say $\hat h=\alpha\hat s_i+\beta\hat s_{i+1}$ with $0\le \alpha,\beta\le D$. We define $B(\hat h)$ as the area of the triangle with sides $\hat h$, $-\alpha\hat s_i$ and $-\beta\hat s_{i+1}$ (in that order). 
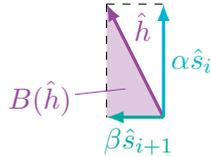
\vspace*{2mm}
\begin{center}
    \begin{tikzpicture}[scale=.25]
	\fill[Purple, opacity=.3] (-9,-9) -- (-6,-9) -- (-9,-3) -- cycle;
		
	\draw[very thick, Purple, -latex] (-6,-9) -- (-9,-3);
    \node[Purple] at (-7.1,-4.2) {$\hat h$};
    
    \draw[very thick, Turquoise, -latex] (-6,-9) -- (-6,-3);
    \draw[dashed] (-9,-3) -- (-6,-3);
    \node[Turquoise] at (-4.5,-6) {$\alpha\hat s_i$};
    
    \draw[very thick, JungleGreen, -latex] (-6,-9) -- (-9,-9);
    \draw[dashed] (-9,-3) -- (-9,-9);
    \node[JungleGreen] at (-7.3,-10.3) {$\beta\hat s_{i+1}$};
	
	\node[Purple, inner sep=1pt] (B) at (-12.5,-8) {$B(\hat h)$};
	\draw[Purple] (B) -- (-8,-7.5);
    \end{tikzpicture} \vspace*{-2mm}
\end{center}
Note that the area $A(h)+B(\hat h)$ is enclosed by a curve of length $\norm{h}_\CC+\|\hat h\|_L\le 2D$ (specifically $h$ concatenated with the two segments $-\alpha\hat s_i$ and $-\beta\hat s_{i+1}$), so that
	\[ \abs{A(h)+B(\hat h)}\le \left(\frac{2D}{\norm z_\CC}\right)^2 \]
by Busemann's isoperimetric inequality.

\newpage

Consider $\gamma_z$ any geodesic representing $z$ and consider the dilation $\delta_\lambda\gamma_z$ with 
\[ \lambda=\sqrt{A(g)+A(h)+B(\hat h)}\ge \sqrt{C\cdot D^4-O(D^2)}=\sqrt{C}\cdot D^2 - O(1). \]
We consider $D\gg 1$ so that all sides of $\delta_\lambda\gamma_z$ have Minkowski-length $\ge D$. In particular we can find two points $h_0$ and $h_1$ on the curve $\delta_\lambda\gamma_z$, more precisely on the sides with direction $-\hat s_{i+1}$ and $-\hat s_i$ respectively, differing by a vector $\hat h$. Up to picking a different starting point on the geodesic $\delta_\lambda\gamma_z$, we may assume $h_0=\mathbf 0$ and $h_1=\hat h$.
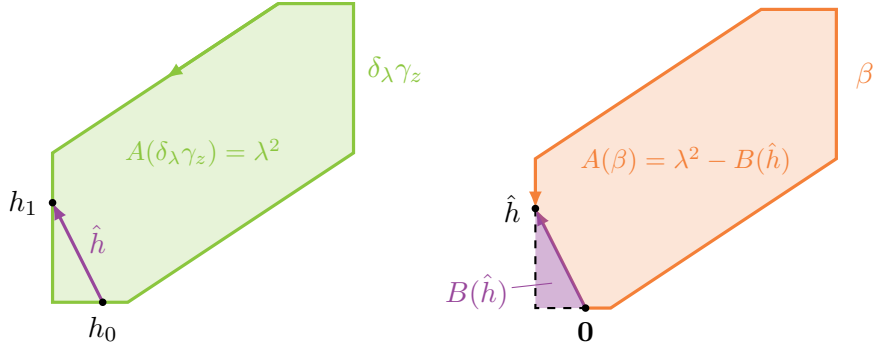
\begin{center}
\begin{tikzpicture}[scale=.22]
	\draw[very thick, LimeGreen, fill=LimeGreen, fill opacity=.2] (9,9) -- (9/2,9) -- (-9,0) -- (-9,-9) -- (-9/2,-9) -- (9,0) -- cycle;
    \draw[very thick, LimeGreen, -latex] (9/2,9) -- (-9/4,9/2);
	\draw[very thick, Purple, -latex] (-6,-9) -- (-9,-3);
 
    \node[circle, fill=black, inner sep=1pt, label=below:$h_0$] at (-6,-9) {};
    \node[circle, fill=black, inner sep=1pt, label=left:$h_1$] at (-9,-3) {};
    
	\node[Purple] at (-6.3,-5.2) {$\hat h$};
    \node[LimeGreen] (B) at (11.5,5) {$\delta_\lambda\gamma_z$};
    \node[LimeGreen] at (0,0) {\small$A(\delta_\lambda\gamma_z)=\lambda^2$};
\end{tikzpicture}
\begin{tikzpicture}[scale=.22]
	\fill[Purple, opacity=.4] (-9,-9) -- (-6,-9) -- (-9,-3) -- cycle;
		
	\draw[thick, dashed] (-9,-3) -- (-9,-9) -- (-6,-9);
	\draw[very thick, Purple, -latex] (-6,-9) -- (-9,-3);
	
    \fill[Orange, fill opacity=.2] (-6,-9) -- (-9/2,-9) -- (9,0) -- (9,9) -- (9/2,9) -- (-9,0) -- (-9,-3);
	\draw[very thick, Orange, -latex] (-6,-9) -- (-9/2,-9) -- (9,0) -- (9,9) -- (9/2,9) -- (-9,0) -- (-9,-3);

    \node[circle, fill=black, inner sep=1pt, label=below:$\mathbf 0$] at (-6,-9) {};
    \node[circle, fill=black, inner sep=1pt, label=left:$\hat h$] at (-9,-3) {};
    
    \node[Purple, inner sep=1pt] (B) at (-12.5,-8) {$B(\hat h)$};
	\draw[Purple] (B) -- (-8,-7.5);
	\node[Orange] at (10.7,5) {$\beta$};
    \node[Orange] at (0,0) {\small$A(\beta)=\lambda^2-B(\hat h)$};
\end{tikzpicture} \vspace*{-2mm}
\end{center}
We consider the curve $\beta$ following $\delta_\lambda\gamma_z$ from $h_0=\mathbf 0$ to $h_1=\hat h$. The curve $\beta$ has endpoint $h_1=\hat h$ (just as $gh$) and area $\lambda^2-B(\hat h)=A(g)+A(h)=A(gh)$, hence represents $gh$.

\bigbreak

It remains to bound the length of $\beta$ to estimate $\norm{gh}_\CC$. We have
	\[ \lambda =\sqrt{n+A(h)+B(\hat h)}\le \sqrt n + \frac12 \frac{A(h)+B(\hat h)}{\sqrt n} \le \sqrt n + \frac{2D^2}{\sqrt n\norm z_\CC^2} \]
and therefore
	\begin{align*}
	\norm{gh}_\CC & \le \norm{\beta}_L = \norm{\delta_\lambda z}_L - \|\hat h\|_L = \lambda\norm z_\CC- \|\hat h\|_L \\
 & \le \norm g_\CC - \|\hat h\|_L + \frac{2D^2}{\sqrt n\norm z_\CC}\le \norm g_\CC - \|\hat h\|_L + \frac2{\sqrt C\norm z_\CC}.
	\end{align*} 
as $\norm g_\CC=\sqrt n \norm{z}_\CC$ and $n\ge C\cdot D^4$. We conclude by setting $C=4\,\big(M^2\norm z_\CC^2\big)^{-1}$.
\end{proof}

\subsection{Almost dead ends in higher Heisenberg groups} \label{sec:higher_heis}
We first recall a construction from \cite{Neumann}, using some formalism from \cite{Stoll}:
\begin{defn}[Centrally amalgamated direct products, see \cite{Stoll}]
	Let $(G_i)_{i=1,\ldots,m},Z$ be groups, and $\iota_i\colon Z\into G_i$ be injective morphisms such that $\iota_i(Z)$ is central in $G_i$. We define
	\[ {\bigotimes}_Z G_i = \faktor{\prod_i G_i\;}{\la \iota_i(z)\iota_j(z)^{-1} \;\big|\; z\in Z \ra}. \]
	Note that each $G_i$ can be naturally identified with a subgroup of ${\bigotimes}_Z G_i$.
\end{defn}
 
 \medbreak
 
 To grasp some intuition, let us talk about $H_{2m+1}(\R) = {\bigotimes}_Z\,G_i$, where $G_i\simeq H_3(\R)$, and $Z=\{z^n\mid n\in\R\}$ is the common center. Once again, elements of $H_{2m+1}(\R)$ can be seen as equivalence classes of curves, more specifically curves in $(\R^2)^m$ starting at $\mathbf 0$. A curve projects on each of the $m$ planes, where we can compute an area. Centrally amalgamating means we add those $m$ areas together, two curves are equivalent if they have same endpoints and sum of areas. This point of view allows to define $\CC$-metrics on $H_{2m+1}(\R)$ for any Minkowski norm on $\R^{2m}$.

 \medbreak
 
\begin{defn} Given metrics $\norm{\,\cdot\,}_i$ on each $G_i$, we can endow $\bigotimes_Z G_i$ with two natural metrics:
\begin{itemize}[leftmargin=8mm]
    \item The \emph{standard-like} metric
    \begin{equation} \label{sec4:length_in_standard}
        \norm{[g_1,\ldots,g_m]}_\stndlike
	 = \min_{\substack{z_1,\ldots,z_m\in Z \\ z_1\ldots z_m=e}}\hspace*{4mm} \sum_{i=1}^m \norm{g_iz_i}_i. \vspace*{-2mm}
    \end{equation}
    \item The \emph{cubical-like} metric
    \begin{equation} \label{sec4:length_in_cubical}
       \norm{[g_1,\ldots,g_m]}_\cubelike
	 = \min_{\substack{z_1,\ldots,z_m\in Z \\ z_1\ldots z_m=e}}\hspace*{2mm} \max_{i=1,\ldots,m} \norm{g_iz_i}_i. 
    \end{equation}
\end{itemize}
\end{defn}
\noindent If we start with word metrics $\norm{\,\cdot\,}_{S_i}$ on each $G_i$, both metrics are word metrics corresponding to
\begin{itemize}[leftmargin=8mm]
    \item the \emph{standard-like} generating set $S_\stndlike={\bigcup}_i\,S_i$.
    \item the \emph{cubical-like} generating set $S_{\cubelike}={\prod}_i(S_i\cup\{e\})$.
\end{itemize}

\bigbreak

\noindent Going back to the $H_{2m+1}(\R)$ example, if we start with $\CC$-metrics $\norm{\,\cdot\,}_{\CC,L_i}$ on each $G_i$, we get two $\CC$-metrics on $H_{2m+1}(\R)$ corresponding to 
\[ L_\stndlike=\ConvHull\Big(\bigcup_i L_i\Big)\subset \R^{2m} \quad\text{and}\quad L_\cubelike=\prod_{i=1}^n L_i \]
respectively. Geodesic curves representing $z^n$ can be understood in both cases:
\begin{itemize}[leftmargin=8mm]
    \item With standard-like metric, any geodesic representing $z^n$ projects to a geodesic for $z^n$ in one of the planes where $\norm z_{\CC,L_i}^2$ is minimal, and to constant curves in all other planes.
    \item With cubical-like metric, any geodesic representing $z^n$ projects to geodesics for some central element $z^{n_i}$ in each plane, and those geodesics have the same length.  
\end{itemize}
Note that, if the $\min$ in (\ref{sec4:length_in_standard}) (resp.\ (\ref{sec4:length_in_cubical})) is achieved for $g_iz_i$'s that are all almost dead ends of depth $D$ and margin $M_i$, then $[g_1,\ldots,g_n]$ is an almost dead end of depth $D$ and margin $\sum_i M_i$ (resp. $\max_i M_i$). This leads to the following generalization of Theorem \ref{almost_dead_ends_in_HR}:

\begin{prop} \label{almost_dead_ends_in_higher}
	Consider the higher Heisenberg group
	\[ H_{2m+1}(\R) = {\bigotimes}_ZG_i \]
	where $G_i\simeq H_3(\R)$, and $Z$ is the common center. Fix polygonal $\CC$-metric $\norm{\,\cdot\,}_{\CC,i}$ on each $G_i$, then Theorem \ref{intro:dead_ends_in_Heis} holds in $(H_{2m+1}(\R),\norm{\,\cdot\,}_{\CC,\cubelike})$, i.e., there exist constants $C,M$ such that, for all $D\gg 1$ and $n\ge C\cdot D^4$, the element $g=z^n$ satisfies
	\[ \forall h\in H_{2m+1}(\R)\text{ such that }\norm h_{\CC,\cubelike}\le D,\quad \norm{gh}_{\CC,\cubelike}\le \norm g_ {\CC,\cubelike}+M. \] 
\end{prop}

\begin{proof} Consider $C_i, M_i>0$ the constants in Theorem \ref{almost_dead_ends_in_HR} for $\big(H_3(\R),\norm{\,\cdot\,}_{\CC,i}\big)$. Let
	\begin{align*}
		C
		& = \max_{i=1,\ldots,m}C_i \cdot \left(\frac{\norm z_{\CC,i}^2}{\norm z_{\CC,1}^2}+\ldots + \frac{\norm z_{\CC,i}^2}{\norm z_{\CC,m}^2}\right), \\
		M
		& = \max_{i=1,\ldots,m} M_i.
	\end{align*}
	Consider $g=z^n$ for $n\ge C\cdot D^4$ and $D\gg 1$. Note that - for $g=z^n$ - formula (\ref{sec4:length_in_cubical}) reads as
	\[ \norm{z^n}_{\CC,\cubelike} = \min_{n_1+\ldots+n_m=n} \; \max_{i=1,\ldots,m}\norm{z^{n_i}}_{\CC,i}. \]
	As $\norm{z^{n_i}}_{\CC,i} = \sqrt{\abs{n_i}}\cdot \norm{z}_{\CC,i}$, the minimum is reached for
	\[ n_i = n \cdot \left(\frac{\norm z_{\CC,i}^2}{\norm z_{\CC,1}^2}+\ldots + \frac{\norm z_{\CC,i}^2}{\norm z_{\CC,m}^2}\right)^{-1} \]
	so that $n_i\ge C_i\cdot D^4$. We can write $h=[h_1,\ldots,h_n]$ for some $h_i\in G_i$ such that $\norm h_{\CC,\cubelike} = \max_{i=1,\ldots,n}\norm{h_i}_{\CC,i}$, and in particular $\norm{h_i}_{\CC,i}\le D$. It follows that
	\[ \norm{gh}_{\CC,\cubelike} \le \max_{i=1,\ldots,n}\norm{z^{n_i}h_i}_{\CC,i}
	 \le \max_{i=1,\ldots,m} \big(\!\norm{z^{n_i}}_{\CC,i} + M_i\big) = \norm{g}_{\CC,\cubelike}+M \]
	 using Theorem \ref{almost_dead_ends_in_HR} for the second inequality.
\end{proof}

\subsection{(Almost) dead ends for word metrics} \label{appendix_dead_end_in_Heis}

In order to get information on $H_{2m+1}(\Z)$ with word metrics, we need some tools to compare the word metric with a well-chosen $\CC$-metric. This materialize in the following theorem

\begin{thm}[Krat \cite{Krat}, Tashiro \cite{Krat_for_non_singular}]
    Let $G$ be a lattice inside a simply connected non-singular nilpotent Lie group $N$. Consider a word metric $\norm{\,\cdot\,}_S$ and the limiting $\CC$-metric $\norm{\,\cdot\,}_\CC$. Then the inclusion $(G,\norm{\,\cdot\,}_S)\into (N,\norm{\,\cdot\,}_\CC)$ is a rough isometry, i.e., there exists $K>0$ s.t.\
    \[ \forall g\in G,\; \norm g_\CC-K \le \norm g_S \le \norm g_\CC+K. \]
\end{thm}
\noindent For instance, we can take $G=H_{2m+1}(\Z)$ inside $N=H_{2m+1}(\R)$. The limiting $\CC$-metric is defined by $L=\mathrm{ConvHull}(\hat S)$, where $\hat S\subset\R^{2m}$ is the projection of $S$ onto the abelianization.

\medbreak

We are now able to prove the following theorem, announced in the introduction:
\begin{thm} \label{almost_dead_ends_in_HZ}
	Consider a pair $(G,S)$ with
 \begin{enumerate}[leftmargin=8mm, label={\normalfont(\alph*)}]
     \item Either $G=H_3(\Z)$ and $S$ any finite symmetric generating set.
     \item Or more generally $G=H_{2m+1}(\Z)$ with $S$ a cubical-like generating set.
 \end{enumerate}
 There exist $C,M>0$ such that, for any $D\gg 1$ and $n\ge C\cdot D^4$, the element $g=z^n$ satisfies
	\[ \forall h\in H_3(\Z) \text{ s.t. }\norm h_S\le D,\quad \norm{gh}_S \le \norm{g}_S + M.\]
\end{thm}
\begin{proof}
	We repeat the proof of Lemma \ref{sec1:rough_dead_ends}. Theorems \ref{almost_dead_ends_in_HR} and \ref{almost_dead_ends_in_higher} hold for specific elements of $H_3(\R)$ and $H_{2m+1}(\R)$, so we get a more precise statement in $H_3(\Z)$ and $H_{2m+1}(\Z)$ too.
\end{proof}

\subsection{A parte - Dead ends in direct products}

\begin{prop}
	Consider $G=\prod_i G_i$ {\normalfont(}i.e., $Z$ trivial above{\normalfont)}, and fix generating sets $S_i$.
	\begin{enumerate}[leftmargin=8mm, label={\normalfont(\alph*)}]
		\item Suppose one of the pairs $(G_j,S_j)$ is infinite and has deep pockets, then so does $(G,S_\cubelike)$.
		
		\item Suppose one of the pairs $(G_j,S_j)$ doesn't have deep pockets, then neither does $(G,S_\stndlike)$.
	\end{enumerate}
\end{prop}
\begin{proof} $\,$\vspace*{2mm}\\
(a) Fix some depth $D$. Let $g_j\in G_j$ be a dead end of depth $\ge D$ and length $\norm{g_j}_{S_j}\ge D$. Consider the element $g=(e,\ldots,e,g_j,e,\ldots,e)\in\prod_i G_i$. For any element $h\in\prod_i G_i$ satisfying
\[ \norm h_\cubelike \le D, \]
i.e., $\norm{h_i}_{S_i}\le D$ for all $i$, we have
\[ \norm{gh}_\cubelike = \max \big\{\! \norm{h_1}_{S_1},\ldots,\norm{g_jh_j}_{S_j},\ldots,\norm{h_m}_{S_m}\!\big\} \le \norm{g_j}_{S_j} = \norm g_\cubelike.   \]
This means that $g$ is a dead end of depth $\ge D$.
\bigbreak
\noindent (b) As $(G_j,S_j)$ doesn't have any deep pocket, there exists some $D$ such that $(G_j,S_j)$ doesn't have any dead end of depth $\ge D$. We show that the same holds in $(\prod_i G_i, S_\stndlike)$. Consider any $g\in\prod_i G_i$. As $g_j\in G_j$ isn't a dead end of depth $D$, there exists $h_j\in G_j$ with $\norm{h_j}_{S_j}\le D$ and $\norm{g_jh_j}_{S_j}>\norm{g_j}_{S_j}$. We define $h=(e,\ldots,e,h_j,e,\ldots,e)\in \prod_i G_i$. We have
\begin{align*}
\norm h_\stndlike & = \,\sum_i \norm{h_i}_{S_i} = 0 + \ldots+ 0 +\norm{h_j}_{S_j}+0+\ldots+0 \le D, \\
\norm{gh}_\stndlike & = \sum_i \norm{g_ih_i}_{S_i} = \sum_i \norm{g_i}_{S_i} +\big(\! \norm{g_jh_j}_{S_j}-\norm{g_j}_{S_j}\!\big) > \norm g_\stndlike.
\end{align*}
so that no $g\in\prod_i G_i$ is a dead of depth $\ge D$.
\end{proof}
\begin{rmk}
Proposition \ref{almost_dead_ends_in_HZ} allows to recover results due to Riley and Warshall in \cite{Riley_Warshall}, namely examples of groups with deep pockets for some generating sets, and without deep pockets for others. Their examples are $(C_2\wr \Z)\times\Z$ and $\Gamma_2(\Z/2\Z)\times \Z$ (here $\Gamma_2(\Z/2\Z)$ is the Baumslag-Remeslennikov group, a finitely presented metabelian group), with what happens to be cubical-like and standard-like generating sets. We propose $H_3(\Z)\times\Z$ (with cubical-like and standard generating sets) as an arguably simpler finitely presented example.
\end{rmk}

	\section{Lamplighters and wreath products}

Let us start with a classical result for lengths in wreath products:
\begin{prop}Consider $G=L\wr Q$, endowed with the standard generating set $S=S_L\cup S_Q$. Elements $g\in L\wr Q$ can be identified with pairs $(\Phi,q)$ with $\Phi\colon Q\to L$ a finitely supported function, and $q\in Q$. Moreover, the length of $g$ is given by the formula
	\[ \norm g_\standard = \sum_{p\in Q} \norm{\Phi(p)}_{S_L} + TS(e_Q;\,\supp(\Phi);\,q), \vspace*{-1mm} \]
	where $TS(x;S;y)$ is the length of the shortest path in the Cayley graph $Cay(Q,S_Q)$ starting at $x$, going through $S$ in some order, and ending at $y$.	
\end{prop}

Using this formula, it was shown in \cite{Dead_ends_Lamplighter} that many of those groups have deep pockets:
\begin{thm}[Theorem 6.1 in \cite{Dead_ends_Lamplighter}] \label{dead_ends_lamplighter}
	Consider $G=L\wr \Z$ with $L$ non-trivial, endowed with the standard generating set $S=S_L\cup\{t^\pm\}$. Suppose that $L$ has dead ends of arbitrary depth w.r.t.\ $S_L$ (for instance if $L$ is finite), then so does $G$ with respect to $S$.
\end{thm}
Note that the condition on $(L,S_L)$ is indeed important, as wreath products like $\Z\wr\Z$ do not have dead ends (at least w.r.t.\ the standard generating set). However, the question of $\NG$-rationality of their complete growth series is still settled by the following proposition:
\begin{thm}\label{saddle_lamplighter}
	Let $G=L\wr\Z$ with $L$ non-trivial. Consider the standard generating set $S=S_L\cup \{t^\pm\}$ with $S_L$ symmetric. There exists a sequence of distinct elements $(g_d)$ satisfying
	\[ \forall h\in G \text{ such that }\norm h_S\le d, \quad \norm{g_dh}_S=\|g_dh^{-1}\|_S. \]
	As a corollary, the associated complete growth series in not $\NG$-rational.
\end{thm}
\begin{proof} Let $\ell \in L\setminus\{e_L\}$ and define
	\[ \Psi_d(q) = \begin{cases} \ell & \text{if } q=\pm d, \\ e_L & \text{otherwise.} \end{cases} \]
	We consider $g_d = (\Psi_d,0)$ i.e., the element with only lamps in a non-trivial state on site $\pm d$,  and the lamplighter guy back at $0$. Consider $h=(\Phi,q)\in L\wr \Z$ with $\norm h_S\le d$. We have $h^{-1}=(\Phi(-q+\,\cdot \,)^{-1},-q)$. Note that $\|h^{-1}\|_S=\norm h_S\le d$ so that
	\[ \supp\Phi\cap \supp\Psi_d = \emptyset = \supp\Phi(-q+\,\cdot\,)^{-1} \cap\supp\Psi_d. \]
	It follows that $\norm{g_dh}_S$ can be easily computed:
	\begin{align*}
		\norm{g_dh}_S
		& = \sum_{p\in \Z} \norm{\Psi_d(p)}_{S_L} + \sum_{p\in \Z}\norm{\Phi(p)}_{S_L} + TS(0;\, B(0,d);\, q) \\
		& = 2 \norm{\ell}_{S_L} + \sum_{p\in \Z}\norm{\Phi(p)}_{S_L} + 2d-\abs q.
	\end{align*}
	(We use that any path going through $\pm d$ must go through the entire interval.) Similarly
	\begin{align*}
		\norm{g_dh^{-1}}_S
		& = \sum_{p\in \Z} \norm{\Psi_d(p)}_{S_L} + \sum_{p\in \Z}\norm{\Phi(-q+p)^{-1}}_{S_L} + TS(0;\, B(0,d);\, -q) \\
		& = 2 \norm{\ell}_{S_L} + \sum_{p\in \Z}\norm{\Phi(p)^{-1}}_{S_L} + 2d-\abs q
	\end{align*}
	Compare both formula, recalling that $S_L$ is symmetric: we have $\norm{g_dh}_S=\norm{g_dh^{-1}}_S$.
\end{proof}
\begin{rmk}
Both Theorem \ref{dead_ends_lamplighter} and Proposition \ref{saddle_lamplighter} can be generalized when $(\Z,\{t^\pm\})$ is replaced by $(Q,S_Q)$ such that the associated Cayley graph is an infinite tree. (So $Q$ is a free product of copies of $\Z$ and $C_2$.) The element $g_d$ can be taken to be $g_d=(\Psi_d,e_Q)$ with
\[ \Psi_d(q) = \begin{cases} \ell & \text{if }\norm q_{S_Q}=d, \\ e_L & \text{otherwise.} \end{cases} \]
However, in these cases, \cite{Parry_Wreath} already concludes non $\Z$-rationality for numerical growth series.
\end{rmk}
	\section{Non algebraicity for $H_3(\Z)$}

\subsection{Context-free languages of polynomial growth and Dyck loops}

A $k$-Dyck word is a word over the set of symbols $[_1,]_1,[_2,]_2,\ldots,[_k,]_k$ satisfying
\begin{itemize}
	\item each symbol appears exactly once,
	\item $[_i$ appears before $]_i$,
	\item and if $[_j$ appears in between $[_i$ and $]_i$, then so does $]_j$.
\end{itemize}
For instance, $[_1[_2]_2[_3[_4]_4]_3]_1[_5]_5$ is a Dyck word while $[_1[_2]_1]_2$ and $]_1$ are not.
\bigbreak
\noindent We fix $\Ac$ an alphabet and a $k$-Dyck word $z$. A $k$-Dyck loop with underlying word $z$ is the set of words obtained by placing fixed words $w_0,\ldots,w_{2k}\in \Ac^*$ in between the parenthesis, and replacing parenthesis $[_i$ and $]_i$ by powers $u_i^{n_i}$ and $v_i^{n_i}$ respectively, with $u_i,v_i\in\Ac^*$ fixed, and $n_i$ any positive integer. For instance,
\[ \big\{ ab(a)^{n_1}bc(ac)^{n_2}ac(c)^{n_2}da(da)^{n_3}(b)^{n_3} abc : n_1,n_2,n_3\in\N \big\} \]
is a $3$-Dyck loop with underlying word $[_1[_2]_2]_1[_3]_3$ obtained with $u_1=a$, $u_2=ac$, $v_2=c$,  $v_1=\varepsilon$, $u_3=da$, $v_3=b$ and $w_0=ab$, $w_1=bc$, $w_2=ac$, $w_3=\varepsilon$, $w_4=da$, $w_5=\varepsilon$, $w_6=abc$.

\begin{defn}
	A language $\Lc\subseteq \Ac^*$ is \emph{$k$-poly-slender} if
	\[ S_\Lc(n) := \#\{w\in\Lc : \ell(w)=n \} = O(n^k). \]
\end{defn}
\noindent In the case of $\CF$ languages, this is equivalent to $B_\Lc(n) := \#\{w\in\Lc : \ell(w)\le n \} = O(n^{k+1})$. We have the following structural result:

\begin{thm}[{Ilie--Rozenberg--Salomaa \cite{poly-slender_CF}}] \label{sec5:Dyck_loops}
A $\CF$ language $\Lc\subseteq \Ac^*$ is $k$-poly-slender if and only if it is a finite union of $(k+1)$-Dyck loops.
\end{thm}
\subsection{Main proof}

\begin{thm} \label{algHeis:Heisenberg NG-alg}
	The complete growth series of $G=H_3(\Z)$ with respect to any finite symmetric generating set $S$ is not $\N G$-algebraic.
\end{thm}
\begin{proof}
	The proof goes by contradiction. We suppose that the complete growth series of $(G,S)$ is $\N G$-algebraic. In particular Theorem \ref{sec5:CF_cross_section} gives a $\CF$ language $\Lc\subset (S\times\N)^*$ evaluating to
	\[ \big\{ (g,\norm g_S) : g\in G \big\}. \]
	Note that the growth $B_\Lc(n)$ is a lower bound on the volume growth of $G$, so that $\Lc$ is poly-slender. It follows from Theorem \ref{sec5:Dyck_loops} that $\Lc$ is a finite union of Dyck loops. Consider one of the Dyck loops forming $\Lc$. When evaluating this language, we get an bunch of elements $(w,\ell)\in G\times \N$ where $\ell$ is the \textit{predicted} length $\norm w_S$. More precisely, we get elements
	\[ g(n_1,n_2,\ldots,n_k) \in G \quad\text{of predicted length}\quad \alpha + n_1\tau_1+ n_2\tau_2 + \ldots + n_k\tau_k \]
	for some fixed $\alpha,\tau_1,\tau_2,\ldots,\tau_k$ (depending only on the chosen Dyck loop). Just to make the notation clear(er), if the word underlying the chosen Dyck loop is ${\color{Red}[_1}{\color{Green}[_2}{\color{blue}[_3]_3}{\color{Purple}[_4]_4}{\color{Green}]_2}{\color{red}]_1}$, then
	\[ g(n_1,n_2,n_3,n_4) = w_0\, {\color{red}u_1^{n_1}}\, w_1\, {\color{Green}u_2^{n_2}}\, w_2\, {\color{blue}u_3^{n_3}}\, w_3\, {\color{blue}v_3^{n_3}}\, w_4\, {\color{Purple}u_4^{n_4}} \, w_5\, {\color{Purple}v_4^{n_4}} \, w_6\, {\color{Green}v_2^{n_2}}\, w_7\, {\color{red}v_1^{n_1}}\, w_8 \]
	for some fixed $u_i,v_i,w_i\in G$. We may always suppose $(u_i,v_i)\ne (e,e)$ (otherwise just forget about those, fuse some $w_j$'s and lower $k$) hence $\tau_i>0$ (otherwise varying $n_i$ gives infinitely many words evaluating in a given sphere in $G$). We may also suppose $w_0=w_1=\ldots=w_{2k-1}=e$ (in which case we will also drop the index for $w_{2k}$). For instance the previous example can be rewritten 
	\[ (\underbrace{w_0\,u_1\,w_0^{-1}}_{\tilde u_1})^{n_1}\;(\underbrace{w_0w_1\,u_2\,(w_0w_1)^{-1}}_{\tilde u_2})^{n_2} \ldots (\underbrace{w_0\ldots w_7\,v_1\, (w_0\ldots w_7)^{-1}}_{\tilde v_1})^{n_1}\cdot (\underbrace{w_0w_1\ldots w_8}_{\tilde w}).\]
	\noindent Note that each element $g\in G$ appears as a $g(n_1,n_2,\ldots,n_k)$ when evaluating one of the Dyck loop. Using the pigeonhole principle we know infinitely many elements of $[G,G]$ appear in a given Dyck loop. From now on we will only consider this specific Dyck loop.
	\medbreak
	\noindent The remainder of the proof goes as follows:
	\begin{enumerate}[leftmargin=8mm, label=(\alph*)]
		\item  We show that $\mathbf 0\in \ConvHull\{ \hat u_i+\hat v_i\}$ (where $\hat g$ is projection of $g$ onto the abelianization).
		
		In particular there exist weights $\lambda_i\ge 0$ such that $\sum_{i=1}^k \lambda_i(\hat u_i+\hat v_i)=\mathbf 0$. As $\hat u_i+\hat v_i\in \Z^2$, we may even suppose $\lambda_i\in\mathbb N$.
		
		\item We treat (i.e., find a contradiction) the case when the underlying word has the form $[_1\cdots]_1$. This is done using \emph{central} almost dead ends of large depth.
		
		\item We treat the other case i.e., when the underlying word contains at least two disjoints pairs of brackets $[_1\cdots]_1[_j\cdots]_j\cdots$. This is done using $\CC$-geometry.
		\begin{enumerate}[leftmargin=6mm, label=(\roman*)]
			\item First we deduce two elements $h_1,h_2\in H_3(\R)\setminus\{e\}$ such that $\hat h_1+\hat h_2=\mathbf 0$ and
			\[ \forall m,n\in\R^+, \quad\norm{\delta_mh_1\cdot \delta_nh_2}_\CC = \norm{\delta_mh_1}_\CC + \norm{\delta_nh_2}_\CC. \]
			
			\item Using Duchin-Mooney's Structure Theorem (Theorem \ref{thm:stru} above), we prove that no such elements exist in $H_3(\R)$.
		\end{enumerate}
	\end{enumerate}
\bigbreak
\noindent (a) Suppose on the contrary $\mathbf 0\notin \ConvHull\{ \hat u_i+\hat v_i\}$, then there exists a linear form $h\colon \R^2\to \R$ such that $h(\hat u_i+\hat v_i)>0$ for all $i$. Let $m>0$ be the minimum of those $k$ values. We have
\[ h\Big( \pi_\ab\big(g(n_1,n_2,\ldots,n_k)\big) \Big) \ge (n_1+n_2+\ldots+n_k)\cdot m + h(\hat w) \]
so that $h(\hat g)>0$ except for finitely many choices of $n_1,n_2,\ldots,n_k$. It follows that $\hat g=\mathbf 0$ (i.e., $g\in [G,G]$) for only finitely many values of the parameters, a contradiction. 
\bigbreak
\noindent (b) Consider $g_0 = g(n_1,n_2,\ldots,n_k)$ any commutator appearing in our Dyck loop. Define
\[ g_n = g(n_1+n\lambda_1,n_2+n\lambda_2,\ldots,n_k+n\lambda_k). \]
Note that all those $g_n$ are distinct commutators. Theorem \ref{almost_dead_ends_in_HZ} tells us that for all $D\gg 1$, there is a number $N_D$ such that for $n\ge N_D$, $g_n$ is an almost dead end of depth $\ge D$ and margin $M$. Now consider for any $n$
\[ \tilde g_n = g(n_1+n\lambda_1+M+1,n_2+n\lambda_2,\ldots,n_k+n\lambda_k).\]
Recall that the underlying word starts and ends by $[_1$ and $]_1$ respectively, and that $g_n$ is a commutator (hence is central in $G$). It follows that we can rewrite
\[ \tilde g_n = u_1^{M+1} \cdot g_n \cdot w^{-1}\,v_1^{M+1}\,w = g_n \cdot u_1^{M+1} \cdot w^{-1}\,v_1^{M+1}\,w\]
Let us take $h = u_1^{M+1}\cdot w^{-1}\,v_1^{M+1}\,w$, $D=\norm h_S$ and $n\ge N_D$. We have $g_n$ and a nearby element $\tilde g_n=g_nh$. We compute their lengths in order to find a contradiction. As both $g_n$ and $\tilde g_n$ appears in the Dyck loop, we know
\begin{align*}
	\norm{g_n}_S & = \alpha + (n_1+n\lambda_1)\tau_1 + \ldots + (n_k+n\lambda_k)\tau_k \\
	\norm{\tilde g_n}_S & = \alpha + (n_1+n\lambda_1+M+1)\tau_1 + \ldots + (n_k+n\lambda_k)\tau_k = \norm{g_n}_S+(M+1)\lambda_1
\end{align*}
We have $\norm{g_nh}_S>\norm{g_n}_S+M$ which is a contradiction with $g_n$ being an almost dead end of depth $D$ and margin $M$.
\bigbreak
\noindent (c-i) Suppose that the underlying word factors into (non-empty) disjoints Dyck words. Write
\[ g(n_1,n_2,\ldots,n_k) = g_1(n_1,n_2,\ldots,n_{j-1}) \cdot g_2(n_j,\ldots, n_k) \cdot w.\]
Fix $m,n\in\N$. We are going to consider $g_1(m\lambda_1,\ldots,m\lambda_{j-1})$ and $g_2(n\lambda_j,\ldots,n\lambda_k)$ and approximate them by better behaved $\delta_m h_1$ and $\delta_n h_2$. In order not to crumble under notations, let us take a specific underlying word $[_1[_2]_2]_1\; [_3]_3[_4]_4$, so that
\[ g_1(n_1,n_2) = u_1^{n_1}u_2^{n_2}v_2^{n_2}v_1^{n_1} \quad\text{and}\quad
g_2(n_3,n_4) = u_3^{n_3}v_3^{n_3}u_4^{n_4}v_4^{n_4} \]
For each $g\in H_3(\R)$, we define $\bar g\in H_3(\R)$ the element represented by a straight line segment from $0$ to $\hat g$. The elements $h_1$ and $h_2$ are defined as follows:
\[ h_1  = \bar u_1^{\lambda_1}\bar u_2^{\lambda_2}\bar v_2^{\lambda_2}\bar v_1^{\lambda_1}
\quad\text{and}\quad
h_2 = \bar u_3^{\lambda_3}\bar v_3^{\lambda_3}\bar u_4^{\lambda_4}\bar v_4^{\lambda_4} \]
\begin{center}
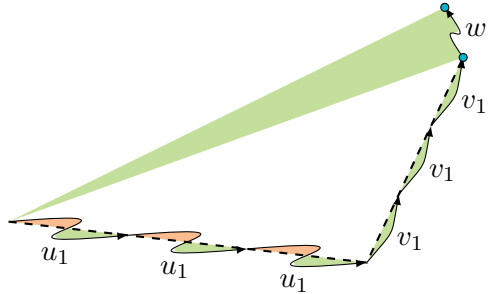

\begin{tikzpicture}[scale=.5, rotate=-25]
\newcommand{\ssqui}[2]{
	\begin{scope}[shift={(#1,#2)}]
	\draw[-latex, looseness=4, out=30, in=210] (0,0) to (3,1);
	\end{scope}}

\newcommand{\tsqui}[2]{
	\begin{scope}[shift={(#1,#2)}]
	\draw[-latex, looseness=1.5, out=75, in=285] (0,0) to (0,2);
	\end{scope}}


\fill[LimeGreen!50] 
(0,0) to[looseness=4, out=30, in=210]
(3,1) to[looseness=4, out=30, in=210]
(6,2) to[looseness=4, out=30, in=210]
(9,3) to[looseness=1.5, out=75, in=285]
(9,5) to[looseness=1.5, out=75, in=285]
(9,7) to[looseness=1.5, out=75, in=285]
(9,9) to (9,3) to (0,0);

\begin{scope}
	\clip (0,0) -- (9,3) -- (9,9);
	\fill[Orange!50] 
		(0,0) to[looseness=4, out=30, in=210]
		(3,1) to[looseness=4, out=30, in=210]
		(6,2) to[looseness=4, out=30, in=210]
		(9,3) to (0,0);
\end{scope}
\draw[thick, dashed] (0,0) -- (9,3) -- (9,9);

\ssqui00
\ssqui31
\ssqui62
\tsqui93
\tsqui95
\tsqui97

\fill[LimeGreen!50] (0,0) to (9,9) to[looseness=3, out=150, in=-30] (8,10) to (0,0);
\draw[-latex, looseness=3, out=150, in=-30] (9,9) to (8,10);

\draw[fill=Turquoise] (9,9) circle (3pt);
\draw[fill=Turquoise] (8,10) circle (3pt);

\node at (1.5,-.3) {$u_1$};
\node at (4.5,.7) {$u_1$};
\node at (7.5,1.7) {$u_1$};
\node at (9.8,4) {$v_1$};
\node at (9.8,6) {$v_1$};
\node at (9.8,8) {$v_1$};
\node at (9,9.8) {$w$};
\end{tikzpicture}
\captionof{figure}{Comparison of $g_1(p)w$ and $\delta_{p}h_1$. (Underlying word starts with $[_1]_1$ and $\lambda_1=1$)}
\end{center}
Let $\pi_\ab\colon H_3(\R)\to \R^2$ be the abelianization map, and recall $A$ is the balayage area. Note that
\begin{align*}
\pi_\ab\big(\delta_{mp}h_1\big) - \pi_\ab\big(g_1(mp\lambda_1,\ldots,mp\lambda_{j-1})\cdot w\big) & = O(1) \\
A\big(\delta_{mp}h_1\big) - A\big(g_1(mp\lambda_1,\ldots,mp\lambda_{j-1})\cdot w\big) & = O(p)
\end{align*}
hence $\norm{\delta_{mp} h_1}_\CC=\norm{g_1(mp\lambda_1,\ldots,mp\lambda_{j-1})\cdot w}_\CC+O(\sqrt{p})$. It follows that
\begin{align*}
	p\norm{\delta_mh_1}_\CC = \norm{\delta_{mp} h_1}_\CC 
	& \sim \norm{g_1(mp\lambda_1,\ldots,mp\lambda_{j-1})\cdot g_2(0,\ldots,0)}_\CC \\
	& \sim \norm{g_1(mp\lambda_1,\ldots,mp\lambda_{j-1})\cdot g_2(0,\ldots,0)}_S \\
	& = \alpha+mp(\lambda_1\tau_1+\ldots+\lambda_{j-1}\tau_{j-1}),
\end{align*}
where the second $\sim$ follows from Krat's theorem, or even Pansu's theorem. Similarly
\begin{align*}
\norm{\delta_{np} h_2}_\CC & =\norm{g_2(np\lambda_j,\ldots,np\lambda_k)\cdot w}_\CC+O(\sqrt{p}) \\
\norm{\delta_{mp} h_1\cdot \delta_{np}h_2}_\CC & =\norm{g_1(mp\lambda_1,\ldots,mp\lambda_{j-1})\cdot g_2(np \lambda_j,\ldots,np\lambda_k)}_\CC+O(\sqrt{p})
\end{align*}
so that
\begin{align*}
p\norm{\delta_nh_2}_\CC & \sim \alpha + np(\lambda_j\tau_j+\ldots+\lambda_k\tau_k) \\
p\norm{\delta_mh_1 \cdot \delta_nh_2}_\CC& \sim \alpha + p\big(m(\lambda_1\tau_1+\ldots+\lambda_{j-1}\tau_{j-1}) + n(\lambda_j\tau_j+\ldots+\lambda_k\tau_k)\big) 
\end{align*}
and finally $\norm{\delta_mh_1\cdot \delta_nh_2}_\CC = \norm{\delta_mh_1}_\CC + \norm{\delta_nh_2}_\CC$ for all $m,n\in\N$ (hence all $m,n\in\R^+$).
\bigbreak
\noindent(c-ii) Suppose that we have two elements $h_1,h_2\ne e$ such that $\hat h_1+\hat h_2=\mathbf 0$ and $\norm{\delta_mh_1\cdot \delta_nh_2}_\CC = \norm{\delta_mh_1}_\CC + \norm{\delta_nh_2}_\CC$ for all $m,n\ge 0$. This means that, for any $\CC$-geodesics $\gamma_1,\gamma_2$ representing $h_1$ and $h_2$, the concatenation $\delta_m\gamma_1\cdot \delta_n\gamma_2$ is a $\CC$-geodesic for $\delta_mh_1\cdot \delta_nh_2$.
\medbreak
\noindent Let us first consider $m=n=1$. We know that $h_1h_2$ is a commutator, so any geodesic has to follow a rescaled isoperimetrix. Suppose wlog that $\norm{h_1}_\CC\ge \norm{h_2}_\CC$, then any geodesic $\gamma_1$ for $h_1$ has to cover at least half the perimeter of the isoperimetrix.
\begin{itemize}[leftmargin=8mm]
	\item If either the isoperimetrix has $\ge 6$ sides or $\norm{h_1}_\CC>\norm{h_2}_\CC$, then $\gamma_1$ has to cover two corners of the isoperimetrix, hence the scale of the isoperimetrix followed by any geodesic continuation of $\gamma_1$ is fixed, and $\gamma_1\gamma_2$ is the maximal geodesic continuation. In particular the longer curve $\gamma_1\cdot\delta_2\gamma_2$ cannot be a geodesic continuation.
	\begin{center}
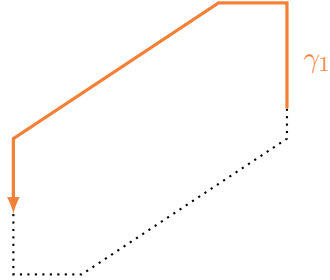

		\begin{tikzpicture}[scale=.2]
		\draw[thick, dotted] (-9,-5) -- (-9,-9) -- (-9/2,-9) -- (9,0) -- (9,2);
		
		\draw[very thick, Orange, -latex] (9,2) -- (9,9) -- (9/2,9) -- (-9,0) -- (-9,-5);
		\node[Orange] at (11,5) {$\gamma_1$};
		\end{tikzpicture}
		\captionof{figure}{All geodesic continuations of $\gamma_1$ have to follow the dotted path.}
	\end{center}

	\item The only remaining case appears when the isoperimetrix has $4$ sides, $\norm{h_1}_\CC=\norm{h_2}_\CC$ and the geodesics $\gamma_1$ and $\gamma_2$ meet at corners of the isoperimetrix. However we still have the same contradiction as $\gamma_1\cdot \delta_2\gamma_2$ is not a $\CC$-geodesic. \qedhere
	\begin{center}
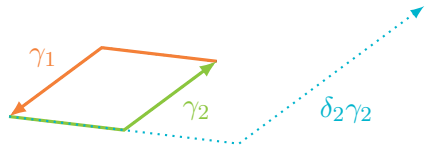

		\begin{tikzpicture}[scale=.2, rotate=-30]
		\draw[thick, dotted] (5,5) -- (-2,2) -- (-5,-5) -- (2,-2) -- (5,5);
		
		\draw[very thick, Orange, -latex] (5,5) -- (-2,2) -- (-5,-5);
		\draw[very thick, LimeGreen, -latex] (-5,-5) -- (2,-2) -- (5,5);
		\draw[thick, Turquoise, dotted, -latex] (-5,-5) -- (9,1) -- (15,15);
		\node[Orange] at (-5,-.5) {$\gamma_1$};
		\node[LimeGreen] at (5.5,1.5) {$\gamma_2$};
		\node[Turquoise] at (14,6.5) {$\delta_2\gamma_2$};
		\end{tikzpicture}
		\captionof{figure}{The path $\gamma_1\cdot\delta_2\gamma_2$ is not quite geodesic.}
	\end{center}
\end{itemize}
\end{proof}
\noindent \textbf{Remark.} This argument extends for $H_{2n+1}(\Z)$ with cubical-like generating sets.
	
	\bibliographystyle{plain}
	\bibliography{bibliography}
\end{document}